\documentclass[Journal]{IEEEtran}

\usepackage{amssymb}
\usepackage{makeidx}         
\usepackage{amsmath}
\usepackage{multicol}        
\usepackage[bottom]{footmisc}
\usepackage{graphicx}
\usepackage{color}
\usepackage{amsfonts}
\usepackage[latin1]{inputenc}
\usepackage{epstopdf}

\newcommand{\eq}{\begin{equation}}
\newcommand{\eeq}{\end{equation}}
\newcommand{\eqn}{\begin{eqnarray}}
\newcommand{\eeqn}{\end{eqnarray}}

\newcommand{\bsea}{\begin{subeqnarray}}
\newcommand{\esea}{\end{subeqnarray}}
\newcommand{\nn}{\nonumber}

\newcommand{\Sp}[2]{\left< #1,#2 \right> }

\newcommand{\Min}[1]{\,\underset{#1}{\mathrm{min}}\,}

\newcommand{\Argmin}[1]{\,\underset{#1}{\mathrm{argmin}}\,}

\newcommand{\Lim}[1]{\,\underset{#1}{\mathrm{lim}}\,}

\newcommand{\de}{\mathrm{d}}
\newcommand{\tr}{\mathop{\rm tr}}  
\newcommand{\e}[1]{\mathrm{e}^{#1}}
\newcommand{\Set}[1]{\left\{ #1\right\}}

\newcommand{\Bc}{ \mathcal{B}}
\newcommand{\Cc}{ \mathcal{C}}
\newcommand{\Dc}{ \mathcal{D}}

\newcommand{\Hc}{ \mathcal{H}}
\newcommand{\Ic}{ \mathcal{I}}

\newcommand{\Kc}{ \mathcal{K}}
\newcommand{\Lc}{ \mathcal{L}}

\newcommand{\Qc}{ \mathcal{Q}}

\newcommand{\Sc}{ \mathcal{S}}

\newcommand{\Vc}{ \mathcal{V}}

\newcommand{\Cs}{ \mathbb{C}}

\newcommand{\Es}{ \mathbb{E}}

\newcommand{\Ns}{ \mathbb{N}}

\newcommand{\Rs}{ \mathbb{R}}
\newcommand{\Ss}{ \mathbb{S}}
\newcommand{\Ts}{ \mathbb{T}}

\newcommand{\Zs}{ \mathbb{Z}}

\def\qed{\hfill \vrule height 7pt width 7pt depth 0pt \smallskip}

\newcounter{pippo}

\newtheorem{remark}{Remark}[section]
\newtheorem{teor}{Theorem}[section]
\newtheorem{corr}{Corollary}[section]
\newtheorem{propo}{Proposition}[section]
\newtheorem{lemm}{Lemma}[section]
\newtheorem{exam}{Example}
\newtheorem{probl}[pippo]{Problem}
\newtheorem{defn}{Definition}[section]
\newcommand{\proof}{\noindent {\bf Proof. }}
\newcommand{\teo}{\begin{teor}}
\newcommand{\eteo}{\end{teor}}
\newcommand{\cor}{\begin{corr}}
\newcommand{\ecor}{\end{corr}}
\newcommand{\prop}{\begin{propo}}
\newcommand{\eprop}{\end{propo}}
\newcommand{\lem}{\begin{lemm}}
\newcommand{\elem}{\end{lemm}}
\newcommand{\ex}{\begin{exam}}
\newcommand{\eex}{\end{exam}}
\newcommand{\pb}{\begin{probl}}
\newcommand{\epb}{\end{probl}}
\newcommand{\df}{\begin{defn}}
\newcommand{\edf}{\end{defn}}
\newcommand{\aprop}{\begin{apropo}}
\newcommand{\eaprop}{\end{apropo}}
\newcommand{\alem}{\begin{alemm}}
\newcommand{\ealem}{\end{alemm}}
\newcommand{\rem}{\begin{remark}}
\newcommand{\erem}{\end{remark}}

\newcommand{\ej}{\mathrm{e}^{j\vartheta}}

\newcommand{\Sts}{\Ss_+^1(\Ts)}
\newcommand{\Stm}{\Ss_+^m(\Ts)}

\newcommand{\Lm}{\mathrm{L}_\infty^{m\times m}(\Ts)}
\newcommand{\Dkl}{\Sc_{\mathrm{KL}}}
\newcommand{\Dis}{\Sc_{\mathrm{IS}}}

\newcommand{\Db}{\Sc_\beta}
\newcommand{\Dnu}{\Sc_{\nu}}
\newcommand{\loga}[1]{\log_{#1}}
\newcommand{\Rgamma}{\mathrm{Range}\;\Gamma}

\newcommand{\Lb}{\Lc_{\nu}}

\newcommand{\Lbg}{\Lc_{\nu}^\Gamma}

\begin{document}

{\color{black}

\title{A New Family of High-Resolution\\ Multivariate Spectral Estimators}

\author{Mattia Zorzi}

\author{Mattia~Zorzi\thanks{This work was supported by the Italian Ministry for Education and
Research (MIUR) under PRIN grant n. 20085FFJ2Z ``New Algorithms
and Applications of System Identification and Adaptive Control"
and by University of Padova under the project ``A Unifying
Framework for Spectral Estimation and Matrix Completion: A New
Paradigm for Identification, Estimation, and Signal Processing".}
\thanks{M. Zorzi is  with the
Department of Electrical Engineering and Computer Science,
University of Liège, Institut Montefiore B28, 4000 Liège, Belgium,
{\tt\small mzorzi@ulg.ac.be}}}

\markboth{DRAFT}{Shell \MakeLowercase{\textit{et al.}}: Bare Demo of IEEEtran.cls for Journals}

\maketitle

\begin{abstract}
In this paper, we extend the Beta divergence family to
multivariate power spectral densities. Similarly to the scalar
case, we show that it smoothly connects the multivariate {\em
Kullback-Leibler} divergence with the multivariate {\em
Itakura-Saito} distance. We successively study a spectrum
approximation problem, based on the Beta divergence family, which
is related to a multivariate extension of the THREE spectral
estimation technique. It is then possible to characterize a family
of solutions to the problem. An upper bound on the complexity of
these solutions  will also be provided. Finally, we will show that
the most suitable solution of this family depends on the specific
features required from the estimation problem.
\end{abstract}


\section{Introduction}
The recent development of THREE-like approaches to multivariate
spectral estimation has triggered a renewed interest for
multivariate distance measures (or simply divergence indexes)
among (power) spectral densities, \cite{DISTANCES_JIANG_2012}. In
the THREE approach, the output covariance of a bank of filters is
used to extract information on the input spectral density. More
precisely, the family of spectral densities matching the output
covariance matrix is considered and a spectrum approximation
problem, which ``chooses'' an estimate of the input spectral
density in this family, is then employed. The choice criterium is
based on finding the spectral density which minimizes a divergence
index with respect to an {\em a priori} spectral density. Note
that, the problem of parameterizing the family of feasible
spectral densities may be viewed as a generalized covariance
extension problem \cite{REALIZATION_GEORGIUO_1987},
\cite{Byrnes95acomplete}, \cite{BYRNES_GEN_ENTROPY_2001}
\cite{GEORGIOU_SPECTRAL_ANALYSIS_2002},
\cite{BYRNES_ONTHEPARTIAL_1997}, \cite{BYRNES_CONVEXOPT_1999}. The
key feature for these estimators concerns the high resolution
achievable in prescribed frequency bands, in particular with short
data records. Significant applications to these methods can be
found in $H_\infty$ robust control
\cite{BYRNES_GENERALIZED_INTERPOLATION_2006},\cite{REMARKS_GEORGIOU_2006},
biomedical engineering \cite{NONINVASIVE_AMINI_2005}, and modeling
and identification \cite{BYRNES_IDENTIFIABILITY_2002},
\cite{GEORGIOU_CONVEX_OPT_2008}, \cite{UNCERTAINTY_KARLSSON_2012}.

The most delicate issue for this theory deals with the choice of
the divergence index. In fact, the corresponding solution to the
spectrum approximation problem (that heavily depends on the
divergence index) must be computable and possibly with bounded
{\em McMillan} degree. Accordingly, it is important to have many
different indexes available in such a way to choose the most
appropriate index in relation to the specific application. The
THREE estimator, introduced by Byrnes, Georgiou and Lindquist in
\cite{A_NEW_APPROACH_BYRNES_2000}, has been extended to the
multichannel case by suggesting different multivariate divergence
indexes, \cite{RELATIVE_ENTROPY_GEORGIOU_2006},
\cite{Hellinger_Ferrante_Pavon},
\cite{FERRANTE_TIME_AND_SPECTRAL_2012}.  In particular, Georgiou
introduced a multivariate version of the {\em Kullback-Leibler}
divergence, \cite{RELATIVE_ENTROPY_GEORGIOU_2006}, which has been
frequently used within information theory, and a multivariate
extension of the {\em Itakura-Saito} distance has been recently
presented by Ferrante {\em et al.},
\cite{FERRANTE_TIME_AND_SPECTRAL_2012}. The latter metric has an
interpretation in terms of relative entropy rate among processes.
Finally, it is worth noting that the output covariance is not
available in a THREE-like spectral estimation method. Indeed, we
need to estimate it by using a collection of sample data generated
by feeding the filters bank with the signal whose spectral density
is to be estimated. Moreover, the family of spectral densities
matching the estimated output covariance must be non-empty. This
covariance estimation task is accomplished by solving a structured
covariance estimation problem,
\cite{ME_ENHANCEMENT_FERRANTE_2012},
\cite{ON_THE_ESTIMATION_ZORZI_2012}. Therefore, a THREE-like
spectral estimation procedure consists in solving a structured
covariance estimation problem and then a spectrum approximation
problem.

The main results of this paper are three. Firstly, we extend to the
multivariate case the Beta divergence family (introduced for the
scalar case in \cite{BASU_ROBUST_1998}) which smoothly connects the
{\em Kullback-Leibler} divergence with the {\em Itakura-Saito}
distance. It is worth mentioning that the Beta divergence family for
scalar spectral densities has been widely used in many applications:
Robust principal component analysis and clustering
\cite{MOLLAH_ROBUST_EXTRACTION_2010}, robust independent component
analysis \cite{MINAMI_ROBUST_BLIND_2002}, and robust nonnegative
matrix and tensor factorization \cite{cichocki_NONNEGATIVE__2009},
\cite{GENRALIZED_ALPHA_BETA_CICHOCKI_2011}.

Secondly, we consider a spectrum approximation problem which
employs the multivariate Beta divergence family. It turns out that
it is possible to characterize a family of solutions to the
problem with bounded {\em McMillan} degree. Moreover, the limit of
the family coincides to the solution obtained by using the {\em
Kullback-Leibler} divergence.

Finally, we tackle the related structured covariance estimation
problem which can be viewed as the static version of the previous
spectrum approximation problem. Also in this case, a Beta matrix
divergence family for covariance matrices, leading to a family of
solutions to the structured covariance estimation problem, may be
introduced.

The paper is outlined as follows. Section \ref{section_THREE_est}
introduces THREE-like spectral estimation methods. Section
\ref{section_beta_family} presents the new extension to the
multivariate case of the Beta divergence family. In Section
\ref{secyion_THREE_beta} the corresponding spectrum approximation
problem is introduced. More precisely, we derive the solution
thanks to the means of the convex optimization. In Section
\ref{section_duale} a non trivial existence result for the dual
problem is established. Then, in Section \ref{section_Newton_alg}
a matricial {\em Newton} algorithm to efficiently solve the dual
problem is presented. In Section \ref{section_sim} some
comparative examples are given: We test the different features of
the found solutions. Section \ref{struct_cov_pb} is devoted to the
estimation of structured covariance matrices by using the Beta
matrix divergence family. Finally, in Section \ref{section_sim2}
we propose an application to the estimation of multivariate
spectral densities which employs the resulting THREE-like
estimator. Moreover, we also draw the different application
scenarios for this family of estimators.

\section{THREE-like spectral estimation} \label{section_THREE_est}
Let us consider an unknown zero mean, $m$-dimensional,
$\Rs^m$-valued, purely non-deterministic, full-rank, stationary
process $y=\Set{y_k;\; k\in\Zs}$ with spectral density
$\Omega(\ej)$ defined on the unit circle $\Ts$. Assume that the
{\em a priori} information on $\Omega$ is given by a {\em prior}
spectral density $\Psi\in\Stm$. Here, $\Stm$ denotes the family of
$\Cs^{m\times m}$-valued spectral density
functions on $\Ts$ which are bounded and coercive, i.e. $\Psi\in\Stm$ if there exist two constants
$\mu_1\geq \mu_2>0$ such that $\mu_2 I\leq\Psi(e^{j\vartheta})\leq \mu_1 I $ on $\Ts$. Then, a finite-length data $y_1 \ldots y_N$
generated by $y$ is observed. We want to find an estimate
$\Phi\in\Stm$ of $\Omega$ by using $\Psi$ and $y_1 \ldots
y_N$. This spectral estimation task is accomplished by employing a
THREE-like approach which hinges on the following four elements:
\begin{enumerate}
  \item A {\em prior} spectral density $\Psi\in\Stm$;
  \item A rational filter to process the data \eq \label{filter_bank}G(z)=(zI-A)^{-1}B,\eeq where
$A\in\Rs^{n\times n}$ is a stability matrix, $B\in\Rs^{n\times m}$
is full rank with $n> m$, and $(A,B)$ is a reachable pair;
  \item An estimate $\hat{\Sigma}$, based on the data $y_1 \ldots y_N$, of the steady
state covariance $\Sigma=\Sigma^T>0$ of the state $x_k$ of the
filter \eq x_{k+1}=Ax_k+By_k;\nn\eeq
  \item A divergence index $\Sc$ between two spectral densities.
\end{enumerate} According to the THREE-like approach, an estimate $\Phi\in\Stm$ of $\Omega$
is given  by solving the problem\footnote{Here and throughout the
paper, integration, when not otherwise specified, is on the unit
circle with respect to the normalized Lebesgue measure. Moreover, a
star denotes transposition plus conjugation.}: \eqn
\label{problema_THREE} && \mathrm{minimize}\; \;\Sc(\Phi\|\Psi)\;\;
\mathrm{over \;the\; set} \nn\\&& \hspace{2cm} \; \;
\Set{\Phi\in\Stm\;|\;\int G\Phi G^*=\hat{\Sigma} }.\eeqn Note that
$\Psi$ is generally not consistent with $\hat{\Sigma}$, i.e. $\int G
\Psi G^*\neq \hat{\Sigma}$. Hence, we have a spectrum approximation
problem. The parametrization of all spectral densities satisfying
constraint in (\ref{problema_THREE}) may be viewed as a generalized
moment problem. For instance, the covariance extension problem may
be recovered by setting \eq \label{filtro_cov_ext_pb}G(z)=\left[%
\begin{array}{ccc}
  z^{-n}I_m & \ldots & z^{-1}I_m \\
\end{array}%
\right]^T.\eeq In this case, the state covariance has a block {\em Toeplitz} structure: \eq\Sigma=\left[%
\begin{array}{ccccc}
  \Sigma_0 & \Sigma_1 &\Sigma_2 &\ldots & \Sigma_{n-1}  \\
  \Sigma_1^T & \Sigma_0 & \Sigma_1 & \ddots & \Sigma_{n-2}\\
  \Sigma_{2}^T & \ddots & \ddots & \ddots & \ddots \\
\Sigma_{n-1}^T & \Sigma_{n-2}^T & \ddots & \ddots & \Sigma_0 \\
\end{array}%
\right],\; \; \Sigma_l=\Es[y_k y_{k+l}^T].\nn\eeq
\subsection{Feasibility of the problem}\label{sezione_feasibility}
The first issue arising with the previous spectrum approximation
problem concerns its feasibility, i.e. the existence of
$\Phi\in\Stm$ satisfying the constraint in (\ref{problema_THREE})
for a given $\hat{\Sigma}$. To deal with this issue, we first
introduce some notation: $\Qc_{n}\subset \Rs^{n\times n}$ denotes
the $n(n+1)/2$-dimensional real vector space of $n$-dimensional
symmetric matrices and $\Qc_{n,+}$ denotes the corresponding cone
of positive definite matrices. We denote as $\Vc(\Ss_+^m)$ the
linear space generated by $\Stm$. Finally, we introduce the linear
operator
\eqn \Gamma &:& \Vc(\Ss_+^m)\rightarrow \Qc_{n} \nn\\
&& \Phi\mapsto \int G\Phi G^*.\nn\eeqn We will see in Section
\ref{section_duale} that the range of $\Gamma$, denoted by
$\Rgamma$, can be profitably exploited for the analysis of the
dual problem of the above spectrum approximation problem. In
\cite{GEORGIUO_THE_STRUCUTRE_2002} (see also
\cite{Hellinger_Ferrante_Pavon}), it was shown that a matrix
$P\in\Qc_{n}$ belongs to $\Rgamma$ if and only if there exists
$H\in \Rs^{m\times
    n}$ such that \eq \label{equazione_R_gamma}P-A P A^T=BH+H^TB^T .\eeq
An equivalent condition, \cite{ME_ENHANCEMENT_FERRANTE_2012}, is
that the kernel of the linear operator \eqn\label{mappa_P}
V&:&\Qc_{n}\rightarrow \Qc_{n}\nn\\ && Q \mapsto
\Pi_B^\bot(Q-AQA^T)\Pi_B^\bot\eeqn contains $P$, namely $V(P)=0$.
Here, $\Pi_B^\bot:=I-B(B^TB)^{-1}B^T$. It turns out that the
spectrum approximation problem is feasible if and only if
$\hat{\Sigma}\in\Rgamma\cap \Qc_{n,+}$,
\cite{GEORGIUO_THE_STRUCUTRE_2002},\cite{Hellinger_Ferrante_Pavon}.
Let $x_1 \ldots x_N$ be the output data generated by feeding the
filters bank with the finite-length data $y_1 \ldots y_N$. An
estimate of $\Sigma$ is therefore given by the sample covariance
matrix $\hat{\Sigma}_C:=\frac{1}{N}\sum_{k=1}^Nx_kx_k^T$ which is
normally positive definite. It may not, however, belong to
$\Rgamma$. Accordingly, we need to compute a new estimate
$\hat{\Sigma}\in\Rgamma$ which is positive definite and ``close''
to the estimate $\hat{\Sigma}_C$. Hence, we have to solve a
structured covariance estimation problem which lead us to consider
the following optimization task. \pb \label{problema_ancillario}
Given $\hat{\Sigma}_C>0$, \eqn && \mathrm{minimize}\;
\;\Dc(P\|\hat{\Sigma}_C)\; \; \mathrm{over\; the
\;set}\nn\\&&\hspace{3cm}\Set{P\in\Qc_{n,+}\;|\; V(P)=0}.\nn\eeqn
\epb Here, $\Dc$ is a suitable divergence index among (positive
definite) covariance matrices. Furthermore, by choosing $\Dc$
convex with respect to $P$, Problem \ref{problema_ancillario} can
be efficiently solved by means of convex optimization. For
instance, in \cite{ME_ENHANCEMENT_FERRANTE_2012} the {\em
information divergence} among two Gaussian densities with
covariance $P$ and $Q$, respectively, \cite{COVER_THOMAS}, has
been considered:
\eq\label{information_divergence}\Dc_\mathrm{I}(P\|Q):=\frac{1}{2}\tr[\log(Q)-\log(P)+PQ^{-1}-I].\eeq
Another approach characterizes $\Sigma$ in terms of the filter
parameters and the sequence of the covariance lags of $y$,
\cite{ON_THE_ESTIMATION_ZORZI_2012}. Once we have $\hat{\Sigma}$
in such a way that the spectrum approximation problem is feasible,
we can replace $G$ with
$\overline{G}=\hat{\Sigma}^{-\frac{1}{2}}G$ and $(A,B)$ with
$(\overline{A}=\hat{\Sigma}^{-\frac{1}{2}}A\hat{\Sigma}^{\frac{1}{2}},\overline{B}=\hat{\Sigma}^{-\frac{1}{2}}B)$.
Thus, the constraint may be rewritten as $\int \bar{G}\Phi
\bar{G}^*=I$. Accordingly, from now on we assume that the spectrum
approximation problem in (\ref{problema_THREE}) is feasible and we
consider the following equivalent formulation. \pb
\label{problema_THREE_I} Given $\Psi\in\Stm$ and
$G(z)=(zI-A)^{-1}B$ such that $I\in\Rgamma$, \eqn
\label{constraint_pb_THREE_I}&&\mathrm{minimize}\; \;
\Sc(\Phi\|\Psi)\; \; \mathrm{over \; the\; set}\nn\\ &&
\hspace{2cm}\; \; \Set{\Phi\in\Stm\;|\; \int G\Phi G^*=I}.\eeqn
\epb

\subsection{Choice of the divergence index}
A divergence index among spectral densities in $\Stm$ must satisfy
the following basic property for all $\Phi,\Psi\in\Stm$:
\eqn\label{proprieta_divergenza} && \Sc(\Phi\|\Psi)\geq 0 \nn\\
&& \Sc(\Phi\|\Psi)=0\; \; \mathrm{if\; and \; only\; if}\; \;
\Phi=\Psi. \eeqn Moreover, the corresponding Problem
\ref{problema_THREE_I} should lead to a computable solution, by
typically solving the dual optimization problem. In
\cite{RELATIVE_ENTROPY_GEORGIOU_2006}, a {\em Kullback-Leibler}
divergence for multivariate spectral densities with the same trace
of the zeroth-moment has been introduced
\eq\label{dist_kl0}\mathcal{S}_{KL0}(\Phi\|\Psi)=\int\tr[\Phi(\log(\Phi)-\log(\Psi))]\eeq
where $\log(\cdot)$, whose definition will be given in Section
\ref{sezione_dist_beta_multi}, is the matrix logarithm. This
divergence is inspired by the {\em Umegaki-von Neumann}'s relative
entropy \cite{NIELSEN_CHUANG_QUANTUM_COMPUTATION} of statistical
quantum mechanics. Moreover, (\ref{dist_kl0}) may be readily
extended to the general case, see
\cite{FLEXIBLE_ROBUST_CICHOKI_AMARI_2010} for the scalar case,
\eq\label{dist_KL}\Dkl(\Phi\|\Psi)=\int
\tr[\Phi(\log(\Phi)-\log(\Psi))-\Phi+\Psi]\eeq and
$\mathcal{S}_{KL0}(\Phi\|\Psi)=\Dkl(\Phi\|\Psi)$ when
$\int\tr\Phi=\int \tr\Psi$. However, the resulting solution to the spectrum
approximation problem is generically non-rational. On the contrary, by considering the multivariate
extension of the {\em Itakura-Saito} distance
\eq\label{D_is}\Dis(\Phi\|\Psi)=\int\tr[
\log(\Psi)-\log(\Phi)+\Phi\Psi^{-1}-I],\nn\eeq the solution is rational when $\Psi$ is rational, \cite{FERRANTE_TIME_AND_SPECTRAL_2012}. We will show in
the following section that the divergence indexes (\ref{dist_KL})
and (\ref{D_is}) belongs to the same multivariate Beta divergence
family. Moreover, this family leads, under a suitable choice of
$\Psi$, to a family of solutions to the spectrum approximation
problem.

Observe that, it is also possible to rewrite Problem
\ref{problema_THREE_I} by considering $\Dkl(\Psi\|\Phi)$. The
resulting solution is, however, only computable when $y$ is a
scalar process
\cite{A_NEW_APPROACH_BYRNES_2000},\cite{KL_APPROX_GEORGIUO_LINDQUIST},
or $\Psi=I$, \cite{GEORGIOU_SPECTRAL_ANALYSIS_2002},
\cite{Blomqvist_MATRIX_VALUED_2003},
\cite{RELATIVE_ENTROPY_GEORGIOU_2006}. Finally, we mention that
there exists another multivariate distance, called {\em Hellinger}
distance, which  gives a rational solution to Problem
\ref{problema_THREE_I}, \cite{Hellinger_Ferrante_Pavon}.

\section{Beta divergence family for spectral densities}\label{section_beta_family}
In this section we extend the notion of Beta divergence (family)
for scalar spectral densities, firstly introduced in
\cite{BASU_ROBUST_1998} and \cite{MINAMI_ROBUST_BLIND_2002}, to
the multivariate case. All the proofs of the propositions stated in this section are placed in Appendix \ref{proofs_SEC_III}.

\subsection{Scalar case}
We recall the definition of the scalar Beta divergence by adopting
the same notation employed in
\cite{FLEXIBLE_ROBUST_CICHOKI_AMARI_2010}. First of all, we need
to introduce the following function  \eqn \loga{c}:& \Rs_+ \times
\Rs_+&\hspace{-0.3cm}\rightarrow \Rs\nn\\&\hspace{-0.7cm}(x,y)&
\hspace{-1cm}\mapsto\left\{
                                                              \begin{array}{ll}
   \frac{1}{1-c}\left(\left(\frac{x}{y}\right)^{1-c}-1\right) , & c\in\Rs\setminus \Set{1} \\
   \log(x)-\log(y) ,&   c=1
                                                              \end{array}
                                                            \right.\nn \eeqn
which is referred to as {\em generalized logarithm discrepancy}
function throughout the paper. Notice that $\loga{c}$ is a
continuous function of real variable $c$ and $\log_c(x,y)=0$ if
and only if $x=y$. The (asymmetric) Beta divergence between two
scalar spectral densities $\Phi,\Psi\in\Sts$ is defined by \eqn
&&\Db(\Phi\|\Psi):=-\frac{1}{\beta}\int\left(\Phi^{\beta}\loga{\frac{1}{\beta}}\left(\Psi^{\beta},\Phi^{\beta}\right)
+\Phi^{\beta}-\Psi^{\beta}\right)\nn\\ &&\hspace{0.5cm}=\int
\left(\frac{1}{\beta-1}\left(\Phi^\beta-\Phi\Psi^{\beta-1}\right)-\frac{1}{\beta}(\Phi^{\beta}-\Psi^\beta)\right)\nn
\eeqn where the parameter $\beta$ is a real number. For $\beta=0$
and $\beta=1$, it is defined by continuity in the following way
\eqn
\lim_{\beta\rightarrow 0}\Db(\Phi\|\Psi)&=&\Dis(\Phi\|\Psi)\nn\\
\lim_{\beta\rightarrow 1}\Db(\Phi\|\Psi)&=&\Dkl(\Phi\|\Psi),\nn\eeqn
where $\Dis$ and $\Dkl$ are the scalar versions of (\ref{D_is})
and (\ref{dist_KL}), respectively. Moreover, the Beta divergence
is a continuous function of real variable $\beta$ in the whole
range including singularities. Thus, it smoothly connects the {\em
Itakura-Saito} distance with the {\em Kullback-Leibler}
divergence. Since property (\ref{proprieta_divergenza}) is
satisfied, $\Db$ is a divergence index. Finally, $\Db$ is always
strictly convex in the first argument, but is often not in the
second argument.
\subsection{Multivariate case} \label{sezione_dist_beta_multi}
Likewise to the scalar case, we start by introducing the {\em
generalized multivariate logarithm discrepancy}. To this aim,
recall that the exponentiation of a positive definite matrix $X$
to an arbitrary real number c, is defined as
$X^c:=U\mathrm{diag}(d_1^c,\ldots,d_m^c)U^T$ where
$X:=U\mathrm{diag}(d_1,\ldots,d_m)U^T$ is the usual spectral
decomposition with $U$ orthogonal, i.e. $UU^T=I$, and
$\mathrm{diag}(d_1,\ldots,d_m)>0$ diagonal matrix.\footnote{It is
also possible to take the exponentiation of positive semidefinite
matrices when $c\neq 0$.} We are now ready to extend the
definition of {\em generalized logarithm discrepancy} to the
multivariate case \eqn \label{def_log_discrepancy_multi}\loga{c}:
&\Qc_{m,+} \times \Qc_{m,+}
&\hspace{-0.3cm}\rightarrow \Rs^{m \times m} \nn\\
&\hspace{-2.5cm}(X,Y) & \hspace{-2cm}\mapsto \left\{
              \begin{array}{ll}
\frac{1}{1-c}\left(X^{1-c}Y^{c-1}-I\right)                , & c\in\Rs\setminus\Set{1} \\
\log(X)-\log(Y), & c=1
    \end{array}
            \right.\eeqn where $\log(X)=U\mathrm{diag}(\log(d_1),\ldots,\log(d_m))U^T$ is the matrix logarithm of
$X$.
 \prop \label{prop_continuity_discrepancy}The {\em generalized multivariate logarithm discrepancy} is a continuous
function of real variable $c$ in the whole range. Moreover,
$\log_{c}(X,Y)=0$ if and only if $X=Y$.\eprop
The exponentiation of a spectral density $\Phi(\ej)\in\Stm$ to an
arbitrary real number $c$ is pointwise defined by using the
previous spectral decomposition:
\eq\label{decomposizione_spettrale_spettri}\Phi(\ej)^c=U(\ej)\mathrm{diag}(d_1(\ej)^c,
\ldots, d_m(\ej)^c) U(\ej)^T\eeq where
$\Phi(\ej)=U(\ej)\mathrm{diag}(d_1(\ej), \ldots, d_m(\ej))
U(\ej)^T$ with $U(\ej)\in\mathrm{L}_\infty^{m\times m}(\Ts)$ such
that $U(\ej)U(\ej)^T=I$. Observe that $\Phi^c$ belongs to $\Stm$.
We are now ready to introduce the multivariate (asymmetric) Beta
divergence among $\Phi,\Psi\in\Stm$: \eqn\label{def_beta_multi}
&& \Db(\Phi\|\Psi):=-\frac{1}{\beta}\int\tr\left[\Phi^{\beta}\loga{\frac{1}{\beta}}\left(\Psi^{\beta},\Phi^{\beta}\right)
+\Phi^{\beta}-\Psi^{\beta}\right]\nn\\&&\hspace{0.5cm}=\int\tr\left[
\frac{1}{\beta-1}(\Phi^\beta-\Phi\Psi^{\beta-1})-\frac{1}{\beta}(\Phi^{\beta}-\Psi^\beta)\right]\eeqn
where $\beta\in\Rs\setminus\Set{0,1}$. Similarly to the scalar
case, we can extend by continuity the definition of Beta
divergence for $\beta=0$ and $\beta= 1$. \prop
\label{prop_limiti_beta_multi}The following limits hold:\eqn
\lim_{\beta\rightarrow
0}\Db(\Phi\|\Psi)&=&\Dis(\Phi\|\Psi)\nn\\
\lim_{\beta\rightarrow 1}\Db(\Phi\|\Psi)&=&\Dkl(\Phi\|\Psi).\nn\eeqn
\eprop
In view of Proposition \ref{prop_continuity_discrepancy} and
Proposition \ref{prop_limiti_beta_multi}, we conclude that the
multivariate Beta divergence is a continuous function of real
variable $\beta$ in the whole range including singularities and it
smoothly connects the multivariate {\em Itakura-Saito} distance
with the multivariate {\em Kullback-Leibler} divergence. \rem For
$\beta=2$, the Beta divergence corresponds, up to a constant
scalar factor, to the standard squared {\em Euclidean} distance
($L_2$-norm) \eq \Sc_{\mathrm{L}2}(\Phi\|\Psi)=\int
\Sp{\Phi-\Psi}{\Phi-\Psi}\nn\eeq where $\Sp{X}{Y}=\tr(XY)$ is the
usual scalar product in $\Qc_m$. \erem Finally, the
multivariate Beta divergence satisfies condition
(\ref{proprieta_divergenza}). \prop \label{proposizione_beta_div_strict_convex} Given $\Phi,\Psi\in\Stm$, the
following facts hold:
\begin{enumerate}
    \item $\Db(\cdot\|\Psi)$ is strictly convex over $\Stm$,
    \item $\Db(\Phi\|\Psi)\geq0$ and equality holds if and
only if $\Psi=\Phi$.
\end{enumerate}       \eprop
Note that $\Db(\Phi\|\cdot)$ is not convex on $\Stm$ (not even in the scalar case).

\section{Spectrum approximation problem}\label{secyion_THREE_beta} Since the Beta divergence is
well-defined for $\beta\in\Rs$, we choose $\beta=-\frac{1}{\nu}+1$
with $\nu\in\Ns_+$\footnote{$\Ns_+$ denotes
the set of the positive natural numbers.}, and we define
$\Dnu(\Phi\|\Psi):=\Sc_\beta(\Phi\|\Psi)$ with
$\beta={-\frac{1}{\nu}+1}$. Moreover, here and in the
remainder part of the paper we assume that
$\Psi(z)^{\frac{1}{\nu}}$ is a rational matrix function. The aim
of this section and Section \ref{section_duale} is to prove
the following statement.
\teo \label{teo_riassuntivo} Given $\Psi\in\Stm$ such that $\Psi(z)^{\frac{1}{\nu}}$ is rational, and $G(z)$ such that $I\in\Rgamma$, the problem
\eq  \label{problem_sec4}\mathrm{minimize}\; \Dnu(\Phi\|\Psi) \;\mathrm{over} \;\Set{\Phi\in\Stm\;|\; \int G\Phi G^*=I}\eeq always admits a unique solution when $\nu\in\Ns_+$. Moreover, such a solution is rational
with {\em McMillan} degree less than or equal to $\nu (\deg[\Psi^{\frac{1}{\nu}}]+2n)$.
\eteo

Since (\ref{problem_sec4}) is a
constrained convex optimization problem,  we consider the
corresponding {\em Lagrange} functional\eqn
&& \hspace{-0.9cm} L_\nu(\Phi, \Lambda)\nn\\ && \hspace{-0.5cm}=\Dnu(\Phi\|\Psi)+\frac{\nu}{1-\nu}\int \tr(\Psi^{\frac{\nu-1}{\nu}}) +\Sp{\int G\Phi G^*-I}{\Lambda}\nn\\
&& \hspace{-0.5cm}= \int \tr\left[-\nu
(\Phi^{\frac{\nu-1}{\nu}}-\Phi\Psi^{-\frac{1}{\nu}})+\frac{\nu}{1-\nu}\Phi^{\frac{\nu-1}{\nu}}
+G^*\Lambda G\Phi\right]\nn\\
&& -\tr[\Lambda]\nn\eeqn
 where we exploited the fact that the term $\int \tr[\Psi^{\frac{\nu-1}{\nu}}]$ plays no role in the optimization problem.
Note that, the {\em Lagrange} multiplier $\Lambda\in\Qc_n$ can be
uniquely decomposed as $\Lambda=\Lambda_\Gamma+\Lambda_\bot$ where
$\Lambda_\Gamma\in\Rgamma$, $\Lambda_\bot\in[\Rgamma]^\bot$. Since
$\Lambda_\bot$ is such that $G^*(\ej)\Lambda_\bot G(\ej)\equiv 0$
and $\tr[\Lambda_\bot]=\Sp{\Lambda_\bot}{I}=0$ (see \cite[Section
III]{MATRICIAL_ALGTHM_RAMPONI_FERRANTE_PAVON_2009}), it does not
affect the {\em Lagrangian}, i.e.
$L_\nu(\Phi,\Lambda)=L_\nu(\Phi,\Lambda_\Gamma)$. Accordingly we
can impose from now on that $\Lambda\in\Rgamma$.

Consider now the unconstrained minimization problem
$\Min{\Phi}\Set{L_\nu(\Phi,\Lambda)\;|\; \Phi\in\Stm}$. Since
$L_\nu(\cdot,\Lambda)$ is strictly convex over $\Stm$, its unique
minimum point $\Phi_\nu$ is given by annihilating its first
variation in each direction $\delta\Phi\in\Lm$:  \eqn &&
\label{variazione_delta_phi_beta}\delta L_\nu(\Phi,\Lambda;\delta
\Phi)\nn\\ && \hspace{0.5cm}=\int\tr\left[\left(\nu(\Psi^{-\frac{1}{\nu}}-\Phi^{-\frac{1}{\nu}})+G^*\Lambda
G\right)\delta\Phi\right] \eeqn where we exploited
(\ref{formula_variazione_esponente}). Note that,
$\nu(\Psi^{-\frac{1}{\nu}}-\Phi^{-\frac{1}{\nu}})+G^*\Lambda G\in
\Lm$. Thus, (\ref{variazione_delta_phi_beta}) is zero $\forall
\delta \Phi\in \Lm$ if and only if \eq
\Phi^{-\frac{1}{\nu}}=\Psi^{-\frac{1}{\nu}}+\frac{1}{\nu} G^*
\Lambda G.\nn\eeq Since $\Phi^{-\frac{1}{\nu}}\in\Stm$, the set of
the admissible {\em Lagrange} multipliers is \eq
\Lb:=\Set{\Lambda\in\Qc_{n}\;|\;
\Psi^{-\frac{1}{\nu}}+\frac{1}{\nu}G^*\Lambda G>0\hbox{ on
}\Ts}.\nn\eeq Therefore, the natural set for $\Lambda$ is \eq
\Lbg=\Lb\cap \Rgamma.\nn\eeq In conclusion, the unique minimum point
of the {\em Lagrange} functional has the form
\eq\label{ottimo_beta}\Phi_\nu(\Lambda):=\left(\Psi^{-\frac{1}{\nu}}+\frac{1}{\nu}G^*\Lambda
G\right)^{-\nu}.\eeq \prop \label{proposizione_grado_mc_millan}
If $\Phi_\nu$ is a minimizer of Problem \ref{problema_THREE_I}, then it is a rational matrix function with {\em McMillan}
degree less than or equal to $\nu (\deg[\Psi^{\frac{1}{\nu}}]+2
n) $. Moreover, the following facts hold:
\begin{enumerate}
    \item If $\Psi$ is constant then, among all the spectral density $\Phi_\nu$ with
$\nu\in\Ns_+$, the spectral density with the
smallest upper bound on the {McMillan} degree corresponds to the
{\em Itakura-Saito} distance
    \item As $\nu\rightarrow +\infty$, $\Phi_\nu$ tends to
    the spectral density corresponding to the {\em Kullback-Leibler} divergence.
\end{enumerate}
\eprop \proof Since both $\Psi^{\frac{1}{\nu}}$ and $G$ are
rational matrix functions and $\nu$ is an integer, then also
$\Phi_\nu$ is rational. Moreover, in view of (\ref{ottimo_beta}),
$\deg[ \Phi_\nu]\leq \nu (\deg[\Psi^{\frac{1}{\nu}}]+2n)$ where
$n$ is the {\em McMillan} degree of $G(z)$.\\1) Since $\Psi$ is
constant, we get $\deg[\Phi_\nu]\leq \nu 2n$ with
$\nu\in\Ns_+$. Thus, the spectral density with the
smallest upper bound on the {\em McMillan} degree is attained for
$\nu=1$, i.e. $\beta=0$, which is the optimal
form related to $\Dis(\Phi\|\Psi)$. Note that,
$\Phi_1(\Lambda)= (\Psi^{-1}+G^*\Lambda G)^{-1}$, which is the
same optimal form found in \cite{FERRANTE_TIME_AND_SPECTRAL_2012}
for the multivariate {\em Itakura-Saito} distance.\\2) Firstly, it is
possible to show that the optimal form obtained by using the {\em
Kullback-Leibler} divergence is
$\Phi_{\mathrm{KL}}(\Lambda)=\e{\log(\Psi)-G^*\Lambda G}$ which is
a straightforward generalization of the optimal form for
$\Sc_{\mathrm{KL}0}$ presented in
\cite{RELATIVE_ENTROPY_GEORGIOU_2006}. We want to show that
$\Phi_\nu\rightarrow \Phi_{\mathrm{KL}}$ as $\nu\rightarrow
+\infty$. Let us consider the function
$F(\lambda):=\log(\Psi^{-\lambda }+\lambda G^*\Lambda G)$ with
$\lambda\in\Rs$ such that $\Psi^{-\lambda }+\lambda G^*\Lambda
G>0$ on $\Ts$. Its first order {\em Taylor} expansion with respect
to $\lambda=0$ is $\Psi^{-\lambda}+\lambda G^*\Lambda G-I$.
Accordingly,
 \eqn && \hspace{-1cm}\Lim{\nu\rightarrow +\infty}
\nu\log\left(\Psi^{-\frac{1}{\nu}}+\frac{1}{\nu}G^*\Lambda G\right)\nn\\ && \hspace{-0.7cm} = \Lim{\nu\rightarrow
+\infty}\frac{\Psi^{-\frac{1}{\nu}}-I}{ \nu^{-1}} +G^*\Lambda
G=-\log(\Psi)+G^*\Lambda G\nn\eeqn where we
exploited (\ref{limite_logaritmo}) and the previous {\em Taylor}
expansion. Finally, \eqn && \Lim{\nu\rightarrow
+\infty}\Phi_\nu(\Lambda)=\Lim{\nu\rightarrow +\infty}
\e{\log\left(\left(\Psi^{-\frac{1}{\nu}}+\frac{1}{\nu}G^*\Lambda
G\right)^{-\nu}\right)}\nn\\ && \hspace{0.5cm}=\Lim{\nu\rightarrow
+\infty} \e{-\nu\log\left(\Psi^{-\frac{1}{\nu}}+\frac{1}{\nu}G^*\Lambda G\right)}\nn\\ && \hspace{0.5cm}=\e{-\Lim{\nu\rightarrow +\infty}
\nu\log\left(\Psi^{-\frac{1}{\nu}}+\frac{1}{\nu}G^*\Lambda
G\right)}\nn\\ && \hspace{0.5cm}=\e{\log(\Psi)-G^*\Lambda
G}=\Phi_{\mathrm{KL}}(\Lambda).\nn\eeqn \qed\\ In the light of
Proposition \ref{proposizione_grado_mc_millan}, there always
exists a unique (up to a right-multiplication by a constant
orthogonal matrix) stable and minimum phase rational spectral
factor $W$ such that $\Psi(\ej)^\frac{1}{\nu}=W(\ej) W(\ej)^*$. By
defining $G_1(\ej)=\frac{1}{\sqrt{\nu}}G(\ej)W(\ej)$, we obtain an
equivalent form of (\ref{ottimo_beta}): \eq
\Phi_\nu(\Lambda)=\left( W(I+G_1^*\Lambda
G_1)^{-1}W^*\right)^{\nu}.\nn\eeq

In this section we showed that $\Phi_\nu(\Lambda)$ is the unique
minimum point of $L_\nu(\cdot,\Lambda)$, namely \eqn &&
\label{disug_lagrangiane}L_\nu(\Phi_\nu(\Lambda),\Lambda)<
L_\nu(\Phi,\Lambda),\; \; \forall \Phi\in\Stm\nn\\ && \hspace{0.5cm}\hbox{s.t.}\;
\Phi\neq \Phi_\nu(\Lambda),\;\Lambda\in\Lbg.\eeqn Hence, if we produce
$\Lambda^\circ\in\Lbg$ such that $\Phi_\nu(\Lambda^\circ)$ is satisfying the constraint in
(\ref{constraint_pb_THREE_I}), inequality
(\ref{disug_lagrangiane}) implies \eq
\Sc_\nu(\Phi_\nu(\Lambda^\circ)\|\Psi)\leq\Sc_\nu(\Phi\|\Psi),\;
\; \forall \Phi\in\Stm \;\hbox{s.t.}\; \int G^*\Phi G=I\nn\eeq and
equality holds if and only if $\Phi= \Phi_\nu(\Lambda^\circ)$,
namely such a $\Phi_\nu(\Lambda^\circ)$ is the unique solution to
Problem \ref{problema_THREE_I} with $\Dnu$. The following step
consists in showing the existence of such a $\Lambda^\circ$ by
using the duality theory.

\section{Dual problem}\label{section_duale}
Here, we deal with the case $\nu\in\Ns_+\setminus\{1\}$, since the existence of
the solution to the dual problem for $\nu=1$ was already showed in
\cite{FERRANTE_TIME_AND_SPECTRAL_2012}. The dual problem
consists in maximizing the functional \eqn \inf_{\Phi}
L_\nu(\Phi,\Lambda)&=&L_\nu(\Phi_\nu,\Lambda)
\nn\\&&\hspace{-2.3cm}=\frac{\nu}{1-\nu}\int\tr\left[\left(\Psi^{-\frac{1}{\nu}}+\frac{1}{\nu}G^*\Lambda
G\right)^{1-\nu}\right]-\tr[\Lambda]\nn\eeqn where we recall that
$\Psi^{\frac{1}{\nu}}$ (and thus also $\Psi^{-\frac{1}{\nu}}$) is
by assumption a rational matrix function. Hence, it is equivalent
to minimize the following functional hereafter referred to as {\em
dual functional}:\eq
J_\nu(\Lambda)=-\frac{\nu}{1-\nu}\int\tr\left[\left(\Psi^{-\frac{1}{\nu}}+\frac{1}{\nu}G^*\Lambda
G\right)^{1-\nu}\right]+\tr[\Lambda].\nn\eeq \teo
\label{teorema_J_stricly_convex}The dual functional $J_\nu$
belongs to $\Cc^\infty(\Lbg)$ and it is strictly convex over
$\Lbg$.\eteo \proof In order to prove
the statement, we need the following first variation of the map
$X\mapsto \tr[X^c]$ (further details may be found in Appendix \ref{exp_matrice}):
\eqn\label{formula_variazione_esponente} &&\delta(\tr[X^c];\delta X)=c\tr[X^{c-1}\delta X].\eeqn  The first variation of
$J_\nu(\Lambda)$ in direction $\delta \Lambda_1 \in\Qc_n$ is \eqn
\label{gradiente_J_beta} &&\hspace{-0.5cm}\delta
J_\nu(\Lambda;\delta\Lambda_1)\nn\\&&
\hspace{0cm}=-\int\tr\left[\left(\Psi^{-\frac{1}{\nu}}+\frac{1}{\nu}G^*\Lambda
G\right)^{-\nu}G^*\delta\Lambda_1 G\right]+\tr[\delta \Lambda_1]\nn\\
&& \hspace{0cm}=-\int\tr\left[\left( W(I+G_1^*\Lambda
G_1)^{-1}W^{*}\right)^{\nu} G^*\delta\Lambda_1 G\right]\nn\\ && \hspace{0.5cm}+\tr[\delta
\Lambda_1].\eeqn The linear form $\nabla J_{\nu,\Lambda}(\cdot):=
\delta J_\nu(\Lambda;\cdot)$ is the {\em gradient} of $J_\nu$ at
$\Lambda$. In order to prove that $J_\nu(\Lambda)\in\Cc^1(\Lbg)$
we have to show that $\delta(J_\nu(\Lambda);\delta\Lambda_1)$, for
any fixed $\delta \Lambda_1$, is continuous in $\Lambda$. To this
aim, consider a sequence $M_n\in\Rgamma$ such that $M_n\rightarrow
0$ and define $Q_N(z)=W(z)(I+G_1(z)^*NG_1(z))^{-1}W(z)^*$ with
$N\in\Qc_n$. By Lemma 5.2 in
\cite{MATRICIAL_ALGTHM_RAMPONI_FERRANTE_PAVON_2009} and since $W$
is bounded on $\Ts$, $Q_{\Lambda+M_n}$ converges uniformly to
$Q_{\Lambda}$. Thus, applying elementwise the bounded convergence
theorem, we obtain \eq \lim_{n\rightarrow \infty } \int G
Q_{\Lambda+M_n}^\nu G^*= \int GQ_{\Lambda}^{\nu} G^*.\nn \eeq
Accordingly,
 $\delta(J_\nu(\Lambda);\delta\Lambda)$ is continuous, i.e. $J_\nu$ belongs to $\Cc^1(\Lbg)$.
In order to compute the second variation, notice that
$Q_\Lambda=\left(\Psi^{-\frac{1}{\nu}}+\frac{1}{\nu}G^*\Lambda
G\right)^{-1}$ and its first variation in direction $\delta
\Lambda\in\Qc_n$ is \eq\delta Q_{\Lambda;\delta
\Lambda}=-\frac{1}{\nu}Q_\Lambda G^*\delta \Lambda G
Q_\Lambda.\nn\eeq Furthermore, consider the operator $\Ic:A\mapsto
A^{\nu}$. By applying the chain rule, we get \eq \delta
(\Ic(A);\delta A)=\sum_{l=1}^{\nu}A^{l-1}\delta A A^{\nu-l}.\nn\eeq
Since \eq \delta J_\nu(\Lambda;\delta \Lambda_1)=-\int
\tr\left[Q^\nu_\Lambda G^*\delta \Lambda_1
G\right]+\tr[\delta\Lambda_1], \nn\eeq the second variation of
$J_\nu$ in direction $\delta\Lambda_1,\delta\Lambda_2\in\Qc_n$ is
\eqn \label{hessiano_J_beta}
&& \delta^2J_\nu(\Lambda;\delta\Lambda_1,\delta\Lambda_2)\nn\\ && \hspace{0.4cm}=-\sum_{l=1}^\nu\int\tr\left[
Q_{\Lambda}^{l-1}\delta Q_{\Lambda;\delta \Lambda_2}
Q_\Lambda^{\nu-l}G^*\delta\Lambda_1 G\right]\nn\\
&&\hspace{0.4cm}=\frac{1}{\nu}\sum_{l=1}^\nu\int\tr\left[
Q_{\Lambda}^{l}G^*\delta\Lambda_2
GQ_\Lambda^{\nu+1-l}G^*\delta\Lambda_1 G\right].\eeqn The bilinear
form $\Hc_{\nu,\Lambda}(\cdot,\cdot)=\delta^2
J_\nu(\Lambda;\cdot,\cdot)$ is the {\em Hessian} of $J_\nu$ at
$\Lambda$. The continuity of $\delta^2 J_\nu$ can be established
by using the previous argumentation. In similar way, we can show
that $J_\nu$ has continuous directional derivatives of any order,
i.e. $J_\nu\in\Cc^k(\Lbg)$ for any $k$. Finally, it remains to be
shown that $J_\nu$ is strictly convex on the open set $\Lbg$.
Since $J_\nu\in\Cc^\infty(\Lbg)$, it is sufficient to show that
$\Hc_\Lambda(\delta \Lambda,\delta\Lambda)\geq 0$ for each $\delta
\Lambda\in\Rgamma$ and equality holds if and only if
$\delta\Lambda=0$. Since $\nu>0$ and the trace of integrands in
(\ref{hessiano_J_beta}) is positive semidefinite when
$\delta\Lambda_1=\delta\Lambda_2$, we have
$\Hc_\Lambda(\delta\Lambda,\delta\Lambda)\geq 0$. If
$\Hc_\Lambda(\delta\Lambda,\delta\Lambda)=0$, then
$G^*\delta\Lambda G\equiv 0$ namely $\delta
\Lambda\in[\Rgamma]^\bot$ (see \cite[Section
III]{MATRICIAL_ALGTHM_RAMPONI_FERRANTE_PAVON_2009}). Since
$\delta\Lambda\in\Rgamma$, it follows that $\delta\Lambda =0$.
 In conclusion,
the Hessian is positive definite and the dual functional is
strictly convex on $\Lbg$. \qed\\

In view of Theorem \ref{teorema_J_stricly_convex}, the dual
problem $\Min{\Lambda}\Set{J_\nu(\Lambda)\;|\; \Lambda\in\Lbg}$
admits at most one solution $\Lambda^\circ$. Since $\Lbg$ is an
open set, such a $\Lambda^\circ$ (if it does exist) annihilates
the first directional derivative (\ref{gradiente_J_beta}) for each
$\delta\Lambda\in\Qc_n$ \eq\Sp{I-\int
G\left(\Psi^{-\frac{1}{\nu}}+\frac{1}{\nu}G^*\Lambda^\circ
G\right)^{-\nu}G^*}{\delta\Lambda}=0\;\; \forall
\delta\Lambda\in\Qc_n\nn \eeq or, equivalently, \eq I=\int
G\left(\Psi^{-\frac{1}{\nu}}+\frac{1}{\nu}G^*\Lambda^\circ
G\right)^{-\nu}G^*=\int G\Phi_\nu(\Lambda^\circ)G^*.\nn\eeq This
means that $\Phi_\nu(\Lambda^\circ)\in\Stm$ satisfies the
constraint in (\ref{constraint_pb_THREE_I}) and
$\Phi_\nu(\Lambda^\circ)$ is therefore the unique solution to
Problem \ref{problema_THREE_I}.

The next step concerns the existence issue for the dual problem.
Although the existence question is quite delicate, since set
$\Lbg$ is open and unbounded, we will show that a $\Lambda^\circ$
minimizing $J_\nu$ over $\Lbg$ does exist. \teo Let
$\nu\in\Ns_+\setminus\Set{1}$, then the dual functional $J_\nu$
has a unique minimum point in $\Lbg$.\eteo \proof Since the
solution of the dual problem (if it does exist) is unique, we only
need to show that $J_\nu$ takes a minimum value on $\Lbg$. First
of all, note that $J_\nu$ is continuous on $\Lbg$, see Theorem
\ref{teorema_J_stricly_convex}. Secondly, we show that
$\tr[\Lambda]$ is bounded from below on $\Lbg$. Since Problem
\ref{problema_THREE_I} is feasible, there exists $\Phi_I\in\Stm$
such that $\int G\Phi_I G^*=I$. Thus, \eqn
\tr[\Lambda]=\tr\left[\int G\Phi_I G^* \Lambda
\right]=\tr\left[\int G^* \Lambda G\Phi_I\right] .  \nn\eeqn Defining
$\alpha=-\nu\tr\int \Psi^{-\frac{1}{\nu}}\Phi_I$, we obtain \eq
\tr[\Lambda]=\nu\tr\left[\int
\left(\Psi^{-\frac{1}{\nu}}+\frac{1}{\nu}G^*\Lambda G\right)\Phi_I
\right]+\alpha.\nn\eeq Since
$\Psi^{-\frac{1}{\nu}}+\frac{1}{\nu}G^*\Lambda G$ is positive
definite on $\Ts$ for $\Lambda\in\Lbg$, there exists a right
spectral factor $\Delta$ such that
$\Psi^{-\frac{1}{\nu}}+\frac{1}{\nu}G^*\Lambda G=\Delta^*\Delta$.
Moreover, $\Phi_I$ is a coercive spectrum, namely there exists a
constant $\mu>0$ such that $\Phi_I(\ej)\geq \mu I$, $\forall \;
\ej\in\Ts$. Starting from the fact that the trace and the integral
are monotonic functions, we get \eqn
\label{traccia_lambda_bounded_below}\tr[\Lambda]&=&\nu\tr\left[\int
\Delta \Phi_I\Delta^*\right]+\alpha\geq\nu \mu\tr\left[\int
\Delta\Delta^*\right]+\alpha\nn\\ &=&\nu \mu\tr\left[\int
\Psi^{-\frac{1}{\nu}}+\frac{1}{\nu}G^*\Lambda
G\right]+\alpha>\alpha\eeqn where we have used the fact that
$\int \tr\left[\Psi^{\frac{1}{\nu}}+\frac{1}{\nu}G^*\Lambda G\right]>0$ when
$\Lambda\in\Lbg$. Finally, notice that
$J_\nu(0)=-\frac{\nu}{1-\nu}\int\tr[\Psi^{\frac{\nu-1}{\nu}}]$.
Accordingly, we can restrict the search of a minimum point to the
set $\Set{\Lambda\in\Lbg\;|\;J_\nu(\Lambda)\leq J_\nu(0)}$. We now
show that this set is compact. Accordingly, the existence of the
solution to the dual problem follows from the Weierstrass'
Theorem. To prove the {\em compactness} of the set, it is
sufficient to show that:
\begin{enumerate}
    \item $\Lim{\Lambda\rightarrow \partial \Lbg} J_\nu(\Lambda)=+\infty $;
    \item $\Lim{\|\Lambda\|\rightarrow \infty} J_\nu(\Lambda)=+\infty
    $.
\end{enumerate}
1) Firstly, recall that $\Psi(z)^{\frac{1}{\nu}}$ is rational by
assumption, thus
$R_\Lambda(z):=\Psi(z)^{-\frac{1}{\nu}}+\frac{1}{\nu}G(z)^*\Lambda
G(z)$ is a rational matrix function. Let $\lambda_{\Lambda,i}(z)$,
$i=1\ldots m$, denote the eigenvalues of $R_\Lambda(z)$. In view
of (\ref{decomposizione_spettrale_spettri}), the eigenvalues of
$R_\Lambda(z)^{1-\nu}$ are $\lambda_{\Lambda,i}(z)^{1-\nu}$.
Moreover, $
\tr(R_\Lambda(z)^{1-\nu})=\sum_{i=1}^m\lambda_{\Lambda,i}(z)^{1-\nu}$
is a rational function because $R_\Lambda(z)$ is a rational matrix
function and $\nu-1\in\Ns_+$. Observe that $\partial\Lbg$ is the
set of $\Lambda\in\Rgamma$ such that $\lambda_{\Lambda,i}(\ej)\geq
0$ on $\Ts$ and there exists $\bar\vartheta$ and $\bar i$ such
that $\lambda_{\Lambda,\bar i}(\e{j\bar\vartheta})=0$.  Thus for
$\Lambda\rightarrow
\partial \Lbg$, $\lambda_{\Lambda,i}(z)^{1-\nu}$ with $i=1\ldots m$ are positive on $\Ts$ and $\lambda_{\Lambda,\bar i}(z)^{1-\nu}$ has
a pole tending to $\e{j\bar\vartheta}\in\Ts$. Accordingly,
$\tr[\int R_\Lambda^{1-\nu}]\geq \int \lambda_{\Lambda,\bar
i}^{1-\nu}\rightarrow \infty$ as $\Lambda\rightarrow
\partial\Lbg$. In view of  (\ref{traccia_lambda_bounded_below}),
we conclude that $J_\nu(\Lambda)=-\frac{\nu}{1-\nu}\tr\left[\int
R_\Lambda^{1-\nu}\right]+\tr[\Lambda]\rightarrow \infty$ as $\Lambda\rightarrow \partial\Lbg$.\\
2) Consider a sequence $\{\Lambda_k\}_{k\in\Ns}\in\Lbg$, such that
\eq \Lim{k\rightarrow \infty}\|\Lambda_k\|=\infty.\nn\eeq Let
$\Lambda^0_k=\frac{\Lambda_k}{\|\Lambda_k\|}$. Since $\Lbg$ is
convex and $0\in\Lbg$, if $\Lambda\in\Lbg$ then
$\xi\Lambda\in\Lbg$ \; $\forall \; \xi\in[0,1]$. Therefore
$\Lambda^0_k\in\Lbg$ for $k$ sufficiently large. Let
$\eta:=\lim\inf \tr[\Lambda_k^0]$. In view of
(\ref{traccia_lambda_bounded_below}),
\eq\tr[\Lambda_k^0]=\frac{1}{\|\Lambda_k\|}\tr[\Lambda_k]>\frac{1}{\|\Lambda_k\|}\alpha\rightarrow
0,\nn\eeq for $k\rightarrow \infty$, so $\eta\geq 0$. Thus, there
exists a subsequence of $\{\Lambda_k^0\}$ such that the limit of
its trace is equal to $\eta$. Moreover, this subsequence remains
on the surface of the unit ball $\partial
\Bc=\Set{\Lambda=\Lambda^T\;|\; \|\Lambda\|=1}$ which is compact.
Accordingly, it has a subsubsequence $\{\Lambda_{k_i}^0\}$
converging in $\partial \Bc$. Let $\Lambda^\infty\in\partial\Bc$
be its limit, thus $\Lim{i\rightarrow\infty}
\tr[\Lambda_{k_i}^0]=\tr[\Lambda^\infty]=\eta$. We now prove that
$\Lambda^\infty\in\Lbg$. First of all, note that $\Lambda^\infty$
is the limit of a sequence in the finite dimensional linear space
$\Rgamma$, hence $\Lambda^\infty\in\Rgamma$. It remains to be
shown that $\Psi^{-\frac{1}{\nu}}+\frac{1}{\nu}G^*\Lambda^\infty
G$ is positive definite on $\Ts$. Consider the unnormalized
sequence $\{\Lambda_{k_i}\}\in\Lbg$: We have that
$\Psi^{-\frac{1}{\nu}}+\frac{1}{\nu}G^*\Lambda_{k_i} G>0$ on $\Ts$
so that
$\frac{1}{\|\Lambda_{k_i}\|}\Psi^{-\frac{1}{\nu}}+\frac{1}{\nu}G^*\Lambda_{k_i}^0
G$ is also positive definite  on $\Ts$ for each $i$. Taking the
limit for $i\rightarrow \infty$, we get that $G^*\Lambda^\infty G$
is positive semidefinite on $\Ts$ so that
$\Psi^{-\frac{1}{\nu}}+\frac{1}{\nu}G^*\Lambda^\infty G>0$ on
$\Ts$. Hence, $\Lambda^\infty\in\Lbg$. Since Problem
\ref{problema_THREE_I} is feasible, there exists $\Phi_I\in\Stm$
such that $I=\int G\Phi_I G^*$, accordingly \eq
\eta=\tr[\Lambda^\infty]=\tr\int G\Phi_I G^*\Lambda^\infty=\tr\int
\Phi_I^{\frac{1}{2}}G^*\Lambda^\infty G\Phi_I^{\frac{1}{2}}.\nn\eeq
Moreover, $G^*\Lambda^\infty G$ is not identically equal to zero.
In fact, if $G^*\Lambda^\infty G \equiv0$, then
$\Lambda^\infty\in[\Rgamma]^\bot$ and $\Lambda^\infty\neq 0$ since
it belongs to the surface of the unit ball. This is a
contradiction because $\Lambda^\infty\in\Rgamma$. Thus,
$G^*\Lambda^\infty G$ is not identically zero and $\eta>0$.
Finally, we have \eqn &&\hspace{-0.7cm} \Lim{k\rightarrow
\infty}J_\nu(\Lambda_k)\nn\\ && =\hspace{-0.2cm}\Lim{k\rightarrow
\infty}-\frac{\nu}{1-\nu}\tr\left[\int
\left(\Psi^{-\frac{1}{\nu}}+\frac{1}{\nu}G^*\Lambda_k
G\right)^{1-\nu}\right]\nn\\ &&\hspace{0.4cm}+\tr[\Lambda_k] \geq
\Lim{k\rightarrow
\infty}\|\Lambda_k\|\tr[\Lambda_k^0]=\eta\Lim{k\rightarrow \infty}
\|\Lambda_k\|=\infty. \nn\eeqn \qed\\ \rem For the case $\nu\in\Zs$ such that $\nu< 0$, (\ref{ottimo_beta})
still holds. Moreover, if (\ref{ottimo_beta}) is a minimizer of Problem \ref{problema_THREE_I}, then it is a rational matrix function with $\deg[\Phi_\nu]\leq |\nu|(\deg[\Psi^{\frac{1}{\nu}}]+2n)$. However, the dual problem may not have solution: The minimum point for
$J_\nu(\Lambda)$ may lie on $\partial\Lbg$, since $J_\nu$ takes
finite values on the boundary of $\Lbg$.\erem


\section{Computation of $\Lambda^\circ$}\label{section_Newton_alg}
We showed that the dual problem always admits a unique solution
$\Lambda^\circ$ on $\Lbg$ for $\nu\in\Ns_+$. In order to find
$\Lambda^\circ$, we use the following matricial {\em Newton}
algorithm with backtracking stage proposed in
\cite{MATRICIAL_ALGTHM_RAMPONI_FERRANTE_PAVON_2009}:
\begin{enumerate}
    \item Set $\Lambda_0=I\in\Lbg$;
    \item At each iteration, compute the {\em Newton} step $\Delta_{\Lambda_i}$
    by solving the linear equation $\Hc_{\nu,\Lambda_i}(\Delta_{\Lambda_i},\cdot)=-\nabla
    J_{\nu,\Lambda_i}(\cdot)$ where, once fixed $\Lambda_i$, $\nabla J_{\nu,\Lambda_i}(\cdot)$
    and $\Hc_{\nu,\Lambda_i}(\cdot,\cdot)$ must be understood as a
    linear and bilinear form of (\ref{gradiente_J_beta}) and (\ref{hessiano_J_beta}), respectively;
    \item Set $t_i^0=1$ and let $t_i^{k+1}=t_i^k/2$ until both of
    the following conditions hold:
    \eqn && \label{cond_appartenenza_set}\hspace{-0.8cm}\Lambda_i+t_i^k\Delta_{\Lambda_i}\in\Lbg\\
    && \label{cond_su_J}\hspace{-0.8cm}J_{\nu}(\Lambda_i+t_i^k\Delta_{\Lambda_i})<J_\nu(\Lambda_i)+\alpha t_i^k\Sp{\nabla J_{\nu,\Lambda_i}}{\Delta_{\Lambda_i}}\eeqn
    with $0<\alpha<1/2$;
    \item Set $\Lambda_{i+1}=\Lambda_i+t_i^k\Delta_{\Lambda_i}$;
    \item Repeat steps 2, 3 and 4 until
    $\|\nabla J_{\nu,\Lambda_i}(\cdot)\|<\varepsilon$ where $\varepsilon$ is a
    tolerance threshold. Then set $\Lambda^\circ=\Lambda_i$.
\end{enumerate}
The computation of the search direction $\Delta_{\Lambda_i}$ is the
most delicate part of the procedure. The corresponding linear
equation reduces to \eq
\label{search_dorection_equation}\frac{1}{\nu}\sum_{l=1}^\nu\int
GQ_{\Lambda_i}^{l}G^*\Delta_{\Lambda_i}
GQ_{\Lambda_i}^{\nu+1-l}G^*=\int GQ_{\Lambda_i}^{\nu}G^*-I\eeq where
$Q_{\Lambda}=W(I+G_1^*\Lambda G_1)^{-1}W^*$. By similar
argumentations used in \cite [Proposition
8.1]{Hellinger_Ferrante_Pavon}, it is possible to prove that there
exists a unique solution $\Delta_{ \Lambda_i}\in\Rgamma$ to
(\ref{search_dorection_equation}). Accordingly, we can easily
compute $\Delta_{\Lambda_i}$ in this way: \begin{enumerate}
    \item Compute \eq \label{integrale_Y}Y=\int GQ_{\Lambda_i}^{\nu}G^*-I;\eeq
\item Compute a
basis $\Set{\Sigma_1 \ldots \Sigma_M}$ for $\Rgamma$ from
(\ref{equazione_R_gamma}) and for each $\Sigma_k$, $k=1\ldots M$,
compute \eq
\label{integrale_Yk}Y_{k}=\frac{1}{\nu}\sum_{l=1}^\nu\int
GQ_{\Lambda_i}^{l}G^*\Sigma_k GQ_{\Lambda_i}^{\nu+1-l}G^*;\eeq
\item Find $\Set{\alpha_k}$ such that $Y=\sum_{k} \alpha_k Y_k$.
Then set  $\Delta_{\Lambda_i}=\sum_{k}\alpha_k \Sigma_k$.
\end{enumerate}
Concerning the evaluation of the integrals in (\ref{cond_su_J}),
(\ref{integrale_Y}) and (\ref{integrale_Yk}), a sensible and
efficient method based on spectral factorization techniques may be
employed. For further details, including the checking of condition
(\ref{cond_appartenenza_set}), we refer to Section VI in
\cite{MATRICIAL_ALGTHM_RAMPONI_FERRANTE_PAVON_2009}.

Finally, it is possible to prove that:
\begin{enumerate}
    \item $J_\nu(\cdot)\in\Cc^\infty(\Lbg)$ is
    strongly convex on the sublevel set $\Kc=\Set{\Lambda\in\Lbg\;|\; J_\nu(\Lambda)\leq J_\nu(\Lambda_0)}$;
    \item The {\em Hessian} is {\em Lipschitz} continuous in
    $\Kc$.
\end{enumerate}
The proof follows the ones in \cite[Section
VII]{MATRICIAL_ALGTHM_RAMPONI_FERRANTE_PAVON_2009} and
\cite[Section VI-C]{FERRANTE_TIME_AND_SPECTRAL_2012} faithfully.
These properties allow us to conclude that the proposed {\em
Newton} algorithm globally converges, \cite[Chapter
9]{BOYD_CONVEX_OPTIMIZATION}. In particular the rate of
convergence is quadratic during the last stage. In this way, the
solution to Problem \ref{problema_THREE_I} may be efficiently
computed.

\section{Simulations results - Part I}\label{section_sim}
In order to test the features of the family of solutions
$\Phi_\nu$ with $\nu\in\Ns_+$, we take into account the following
comparison procedure:
\begin{enumerate}
    \item Choose a zero mean stationary process $y=\{y_k;k\in\Zs\}$
    with spectral density $\Omega\in\Stm$;
    \item Design a filters bank $G(z)$ as in (\ref{filter_bank});
    \item Choose a prior spectral density $\Psi\in\Stm$ such that
    $\Psi^{\frac{1}{\nu}}$ is rational
    \item Set $\hat{\Sigma}=\Sigma\in \Rgamma\cap\Qc_{n,+}$,
    i.e. the corresponding spectrum approximation problem is feasible;
    \item Solve Problem \ref{problema_THREE_I} (with $\Dnu$) by means of the
    proposed algorithm with the chosen $\Psi$ and $\hat{\Sigma}^{-\frac{1}{2}}G(z)$ as filters bank.
 \end{enumerate}
In the above comparison procedure we assume to know $\Sigma$. In
this way, we avoid the approximation errors introduced by the
estimation of $\Sigma$ from the finite-length data $y_1 \ldots
y_N$. As noticed in Section \ref{section_THREE_est}, $\Psi$
incorporates the {\em a priori} information on $y$. More
specifically, $\Psi$ is designed by using some given partial
information on $y$ (e.g. its zeroth moment), or given laws (e.g.
physical laws if $y$ describes a physical phenomenon) which
describe its theoretical features. When no {\em a priori}
information is available, we choose $\Psi= I$ which represents the
spectral density of the most unpredictable random process.
Concerning the design of the filter, its role consists in
providing the interpolation conditions for the solution to the
spectrum approximation problem. More specifically, a higher
resolution can be attained  by selecting poles in the proximity of
the unit circle, with arguments in the range of frequency of
interest, \cite{A_NEW_APPROACH_BYRNES_2000}.

\subsection{Scalar case}\label{prima_sez_simu_scal}
\begin{figure*}[htbp]
\begin{center}
\includegraphics[width=15cm]{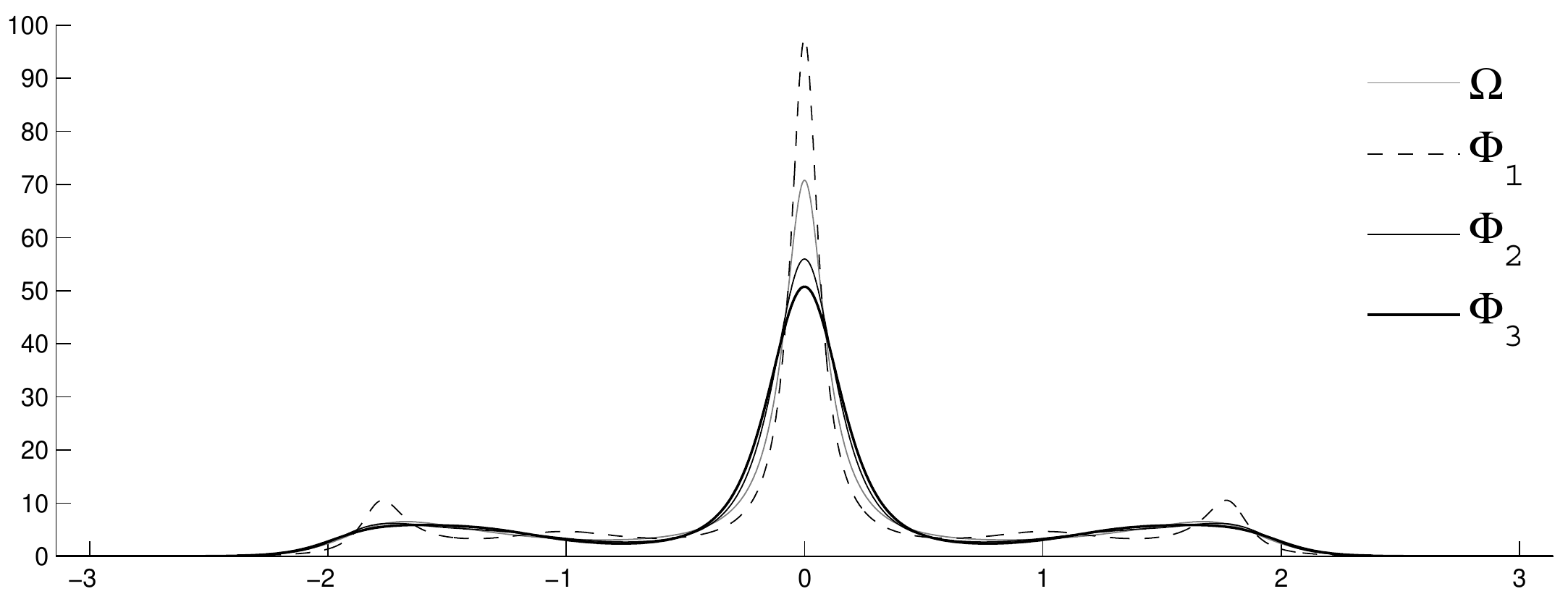}
\end{center}
\vspace{-6mm} \caption{Approximation of an ARMA $(6,4)$ spectral
density.}\label{foto3}
\end{figure*}
We start by considering Example described in \cite[Section
VIII-B]{MATRICIAL_ALGTHM_RAMPONI_FERRANTE_PAVON_2009} (the unique
difference is that we assume to know $\Sigma$ and $\int\Omega$).
Consider the following ARMA process: \eqn
y(t)&=&0.5y(t-1)-0.42y(t-2)+0.602y(t-3)\nn\\
&-& 0.0425 y(t-4)+0.1192y(t-5)\nn\\
&+& e(t) +1.1e(t-1)+0.08e(t-2)\nn\\ &-&0.15e(t-3)\nn \eeqn where $e$ is a
zero-mean {\em Gaussian} white noise with unit variance. In Figure
\ref{foto3}, the spectral density $\Omega\in\Sts$ of the ARMA
process is depicted (gray line). $\Psi$ is equal to $\int \Omega$
and $G(z)$ is structured according to the covariance extension
setting (\ref{filtro_cov_ext_pb}) with $6$ covariance lags (i.e.
$n=6$). In Figure \ref{foto3} the different solutions obtained by
fixing $\nu=1$, dashed line, $\nu=2$, solid line, and $\nu=3$,
thick line, are shown. The solution obtained by minimizing the
multivariate {\em Itakura-Saito} distance ($\nu=1$) is
characterized by peaks which are taller than these in $\Omega$.
On the contrary, the peaks are reduced by increasing $\nu$.
Finally, the solutions with $\nu=2$ and $\nu=3$ are closer to
$\Omega$ than the one with $\nu=1$.

As second example we consider the scalar {\em bandpass} random
process with spectral density $\Omega$ depicted in Figure
\ref{foto1}
\begin{figure}[htbp]
\begin{center}
\includegraphics[width=8cm]{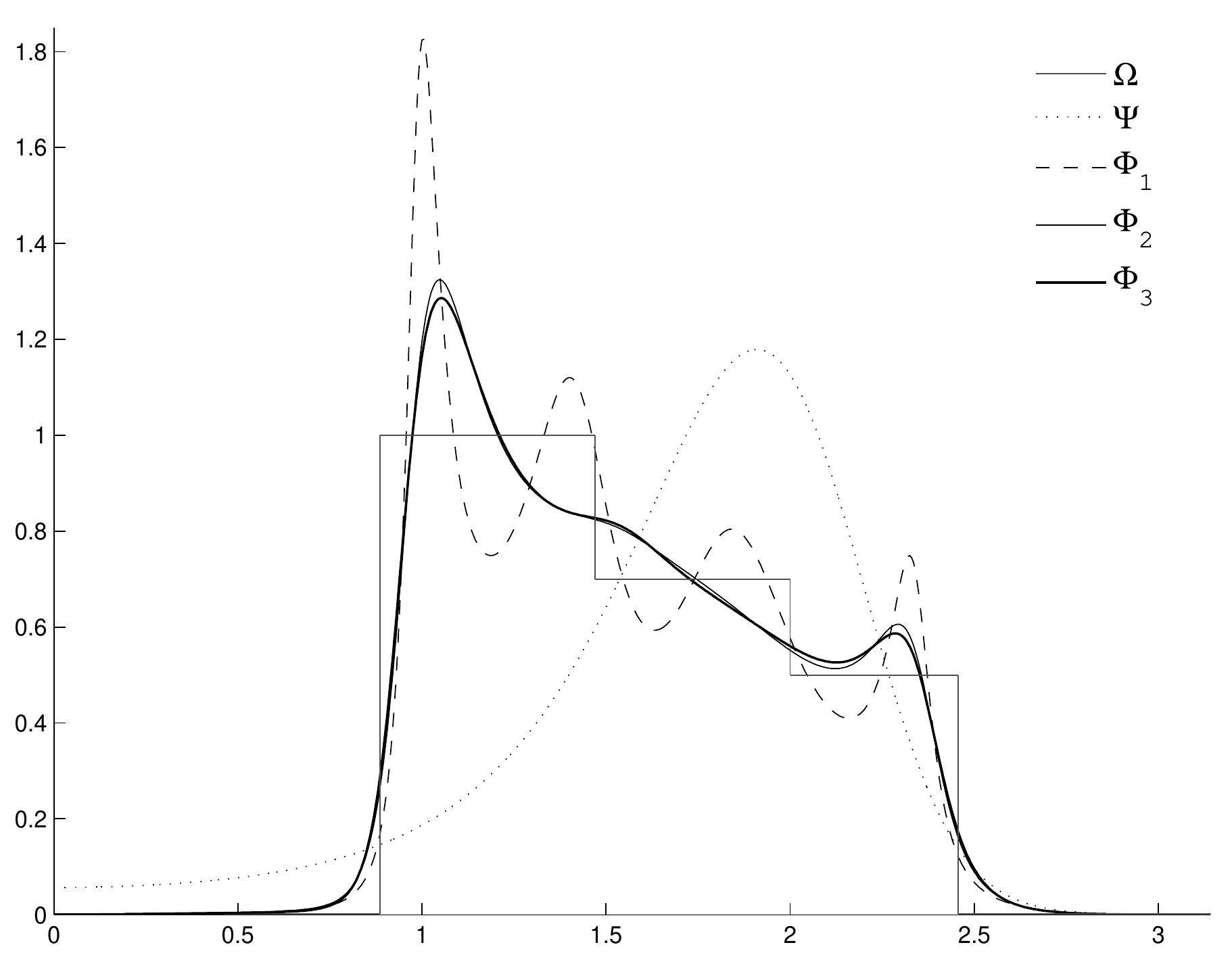}
\end{center}
\vspace{-6mm} \caption{Approximation of the spectral density of a
scalar {\em bandpass} random process.}\label{foto1}
\end{figure} (gray curve). The {\em cutoff} frequencies are $\vartheta_1=0.89$
and $\vartheta_2=2.46$. Moreover, $\Omega(\ej)\geq 2\cdot10^{-3} $
in the {\em stopband}, accordingly $\Omega\in\Sts$. Matrix $B$ is
a column of ones. Matrix $A$ is chosen as a block-diagonal matrix
with one eigenvalue equal to zero and eight eigenvalues equispaced
on the circle of radius $0.8$
 \eq\pm 0.8,\; 0.8\e{\pm j\frac{\pi}{4}},\; 0.8\e{\pm j \frac{\pi}{2}},\; 0.8\e{\pm j\frac{3}{4}\pi
}.\nn\eeq Here, \eqn && \Psi(z)=(W_\Psi(z)W_\Psi(z^{-1}))^6,\nn\\
&& W_\Psi(z)=\frac{5}{6}\frac{z+0.6}{(z-0.4\e{j2.3})(z-0.4\e{-j2.3})}.
\nn\eeqn In this way $\Psi^{\frac{1}{\nu}}$ with $\nu=1$, $\nu=2$, and
$\nu=3$ are rational. Figure \ref{foto1} also shows $\Psi$ and the
obtained solutions. The one with $\nu=1$ turns out inadequate. The
solutions with $\nu=2$ and $\nu=3$ are, instead, similar and
closer to $\Omega$.

\subsection{Multivariate case}\label{subsection_simu_multi}
We consider a bivariate {\em bandpass} random process with
spectral density $\Omega$ plotted in Figure \ref{foto2} (gray
curve).
\begin{figure}[htbp]
\centering
\begin{center}
\includegraphics[width=9.5cm]{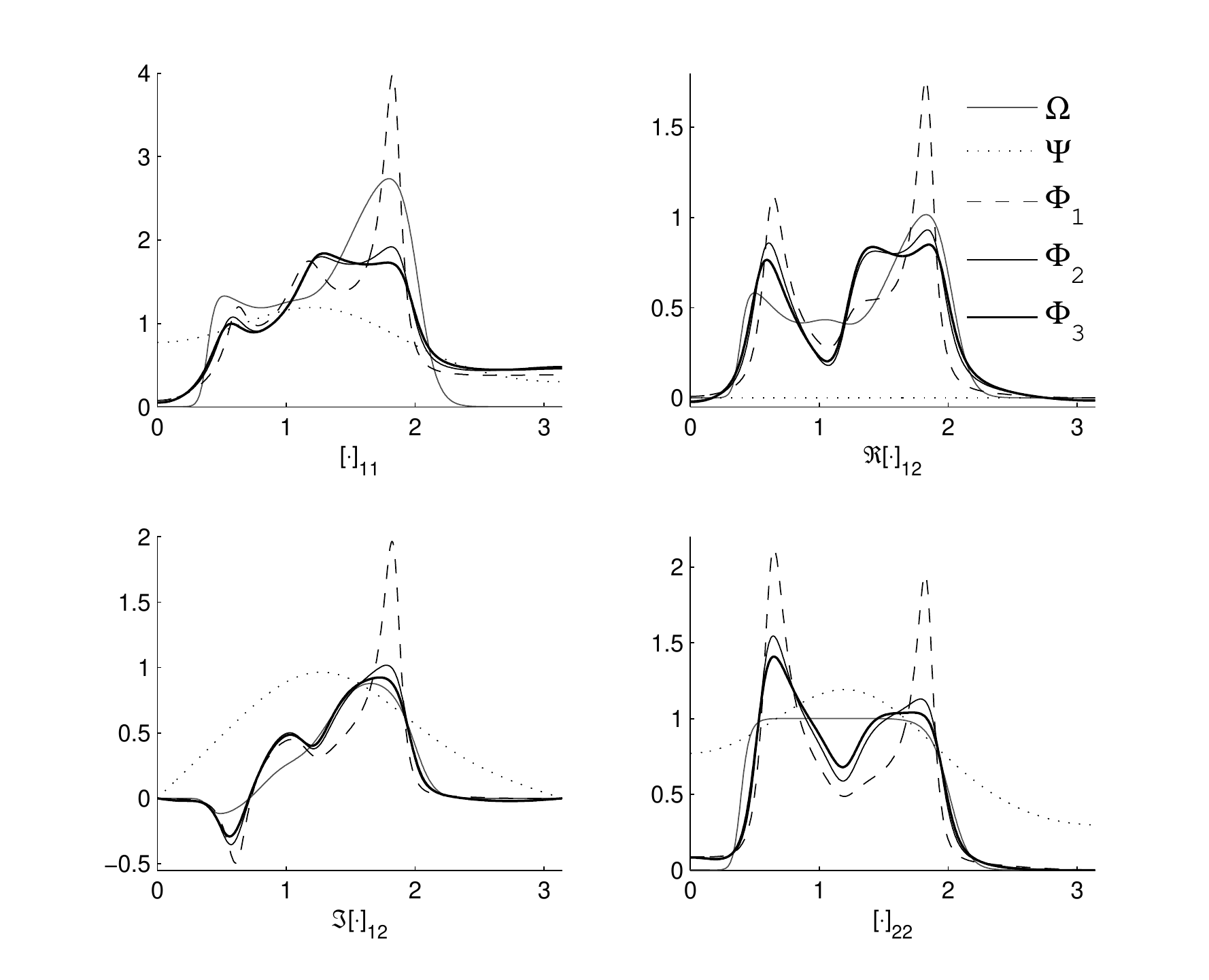}
\end{center}
\vspace{-6mm} \caption{Approximation of the spectral density of a
bivariate {\em bandpass} random process.}\label{foto2}
\end{figure}
Here, the {\em cutoff} frequencies are $\vartheta_1=0.42$ and
$\vartheta_2=1.94$, and $\Omega(\ej)\geq 2\cdot10^{-3} I$ in the
whole range of frequencies. The {\em prior} $\Psi$ is depicted in
Figure \ref{foto2} (dotted line). The matrix $A$ of the filters
bank has one eigenvalue equal to zero, two eigenvalues in $\pm
0.8$ and three pairs of complex eigenvalues closer to the passband
$0.8\e{\pm j0.4},\; 0.8\e{\pm j 1.2},\; 0.8\e{\pm j 2}$. The
solutions for $\nu=1$ (dashed line) $\nu=2$ (solid line) and
$\nu=3$ (thick line) are shown in Figure \ref{foto2}. It is
apparent that the solutions for $\nu=2$ and $\nu=3$ are the most
appropriate.


In view of the previous examples, we now try to point out the
features of the family of solutions. In the above examples the
chosen priors are not characterized by peaks. The found solutions,
however, exhibit peaks which are reduced by increasing $\nu$. In
order to give an interpretation of such a result, consider two
scalar spectral densities $\Psi,\Phi\in\Sts$. Let
$\psi=\Psi(\e{j\bar\vartheta})$, $\phi=\Phi(\e{j\bar \vartheta})$,
where $\bar\vartheta\in [0,2\pi)$ is fixed, and consider the
following function \eqn
&& \hspace{-0.5cm}s_\nu(\phi,\psi)\nn\\ && \hspace{-0.5cm}=\left\{%
\begin{array}{ll}
\hbox{\small $\log \psi-\log \phi+\phi\psi^{-1}-1,$} & \hbox{{\small $\nu=1$}}. \\
  \hbox{\small $-\nu(\phi^{\frac{\nu-1}{\nu}}-\phi\psi^{-\frac{1}{\nu}})-\frac{\nu}{\nu-1}(\phi^{\frac{\nu-1}{\nu}}-\psi^{\frac{\nu-1}{\nu}})$,} & \hbox{{\small $1<\nu<\infty$}} \\
\hbox{\small $\phi(\log\phi-\log\psi)-\phi+\psi,$} & \hbox{{\small $\nu=\infty$}}\\
\end{array}%
\right.\nn \eeqn Informally stated, $s_\nu$ represents the
(infinitesimal) contribution at $\bar\vartheta$ to
$\Sc_\nu(\Phi\|\Psi)$. Note that $s_\nu(\psi,\psi)=0$ for each $\nu$. Since $\Psi$ is given in Problem \ref{problema_THREE_I}, we assume that $\psi$
is a fixed parameter and we consider \eqn && s^\prime_\nu(\phi,\psi):=\left. \frac{\mathrm{d}s_\nu(x,\psi)}{\mathrm{d}x}\right|_{x=\phi}\nn\\ && \hspace{0.5cm}=
 \left\{
   \begin{array}{ll}
     \nu(\psi^{-\frac{1}{\nu}}-\phi^{-\frac{1}{\nu}}), & 1\leq \nu<\infty \\
     \log\phi-\log\psi, & \nu=\infty
   \end{array}
 \right.\nn
\eeqn which represents the instantaneous rate of change of $s_\nu(\cdot,\psi)$ at point $\phi$. The first {\em Taylor} expansion of $s_\nu^\prime(\cdot,\psi)$
with respect to $\phi=\psi$ is the straight line \eq \left\{
                                                       \begin{array}{ll}
                                                         \psi^{-1-\frac{1}{\nu}}(\phi-\psi), & 1\leq    \nu<\infty \\
                                                         \psi^{-1}(\phi-\psi), & \nu=\infty
                                                       \end{array}
                                                     \right.\nn
\eeq
having a slope equal to $\psi^{-1-\frac{1}{\nu}}$ when $1\leq \nu<\infty$ and $\psi^{-1}$ for $\nu=\infty$. Once $\nu$ is fixed, the slope decreases as $\psi$ increases, and it is close to zero for $\psi$ sufficiently large. Thus the critical cases, i.e. when $s_\nu(\phi,\psi)$ is not able to discriminate $\phi$ from $\psi$ sufficiently well, happen when $\psi$ is large,  because $s_\nu(\cdot,\psi)$ is almost flat in a neighborhood of $\psi$.
On the other hand, if $\psi$ is greater than one then the slope increases as $\nu$ increases, i.e. $s_\nu(\phi,\psi)$ improves the ability to discriminate $\phi$ from $\psi$ by increasing $\nu$. Accordingly, a sufficiently large value of $\nu$ avoids solutions
$\Phi$ which are very different from $\Psi$ in narrow ranges of
frequencies. This explains the presence of relevant peaks only for $\nu=1$. The same
conclusion can be obtained by considering multivariate spectral
densities. Concerning the complexity upper bound of the found
solutions, in view of Proposition
\ref{proposizione_grado_mc_millan}, it is easy to check that the
upper bound on the {\em McMillan} degree of $\Phi_\nu^\circ$
increases as $\nu$ increases. For instance, in the first example
of Section \ref{prima_sez_simu_scal} we have
$\deg[\Phi_1^\circ]\leq 12$, $\deg[\Phi_2^\circ]\leq 24$,
$\deg[\Phi_3^\circ]\leq 36$. Thus, the solution with $\nu=1$
guarantees a simple model for the process $y$. Finally, we require
that $\Psi^{\frac{1}{\nu}}$ is rational, accordingly the solution
with $\nu=1$ is the most appropriate to incorporate rational
priors.

\section{Structured covariance estimation problem}\label{struct_cov_pb}
As mentioned in Section \ref{sezione_feasibility}, we only have a
{\em prior} $\Psi$ and a finite-length data $y_1 \ldots y_N$ in the
THREE-like spectral estimation procedure. Moreover, $\Phi_\nu$
represents a family of estimates of $\Omega$ and we showed how to
compute it starting from $\Psi$ and
$\hat{\Sigma}\in\Rgamma\cap\Qc_{n,+}$. Accordingly, it remains to
find $\hat{\Sigma}$ from $y_1 \ldots y_N$. To deal with this issue,
we consider Problem \ref{problema_ancillario} which can be viewed as
the static version of Problem \ref{problema_THREE_I}. Indeed, in
both problems minimization of a divergence index, with respect to
the first argument, is performed on the intersection among a vector
space and an open cone. In this section, we briefly show that it is
also possible to find a family of solutions to the structured
covariance estimation problem.

The Beta matrix divergence (family) among two covariance matrices
$P,Q\in\Qc_{n,+}$ with $\beta\in\Rs\setminus \Set{0,1}$ is defined
as \eq
\Dc_{\beta}(P\|Q):=\tr\left[\frac{1}{\beta-1}(P^\beta-PQ^{\beta-1})-\frac{1}{\beta}(P^\beta-Q^\beta)\right].\nn\eeq
In fact, $\Dc_\beta(P\|Q)$ is the Beta divergence
$\Sc_\beta(\Phi\|\Psi)$ among the two constant spectral densities
$\Phi(\ej)\equiv P$ and $\Psi(\ej)\equiv Q$. Since $\Dc_\beta$ is
a special case of $\Sc_\beta$, it is strictly convex with respect
to the first argument. Moreover, it is a continuous function of
real variable $\beta\in\Rs$ with \eqn
\lim_{\beta\rightarrow 0}\Dc_\beta(P\|Q)&=&\Dc_\mathrm{B}(P\|Q)\nn\\
\lim_{\beta\rightarrow
1}\Dc_\beta(P\|Q)&=&\Dc_{\mathrm{KL}}(P\|Q)\nn\eeqn where
$\Dc_{\mathrm{B}}=2\Dc_{\mathrm{I}}$ (see
(\ref{information_divergence})) is the {\em Burg matrix }
divergence, and \eq
\Dc_{\mathrm{KL}}(P\|Q):=\tr\left[P(\log(P)-\log(Q))-P+Q\right]\nn\eeq
is the extension of the {\em Umegaki-von Neumann}'s relative
entropy, \cite{NIELSEN_CHUANG_QUANTUM_COMPUTATION}, to non
equal-trace matrices.

Take into account Problem \ref{problema_ancillario} with
$\Dc_\nu(P\|\hat{\Sigma}_C):=\Dc_\beta(P\|\hat{\Sigma}_C)$ such
that $\beta=-\frac{1}{\nu}+1$ and $\nu\in\Ns_+$. In
\cite{ME_ENHANCEMENT_FERRANTE_2012}, the existence and uniqueness
of the solution to the problem for $\nu=1$ has been showed.
Moreover, the form of the optimal solution is
$P_{\mathrm{B}}(\Delta)=\left(\hat{\Sigma}_C^{-1}+V^\star(\Delta)\right)^{-1}$,
where $V^\star(\Delta):=\Pi_B^\bot\Delta\Pi_B^\bot-A^T\Pi_B^\bot
\Delta \Pi_B^\bot A$ is the adjoint operator of the linear map $V$
defined in (\ref{mappa_P}) and $\Delta\in\Qc_n$ is the {\em
Lagrange} multiplier. Consider now Problem
\ref{problema_ancillario} with $\nu\in\Ns_+\setminus\Set{1}$. The
corresponding {\em Lagrange} functional is \eqn
L_\nu(P,\Delta)&:=&\Dc_\nu(P\|\hat{\Sigma}_C)+\frac{\nu}{1-\nu}\tr\left[\hat{\Sigma}_C^{\frac{\nu-1}{\nu}}\right]+\Sp{V(P)}{\Delta}\nn\\
&=&\Dc_\nu(P\|\hat{\Sigma}_C)+\frac{\nu}{1-\nu}\tr\left[\hat{\Sigma}_C^{\frac{\nu-1}{\nu}}\right]\nn\\ && \hspace{0.2cm}+\Sp{P}{V^\star(\Delta)}.\nn\eeqn
Since $L_\nu(P,\Delta+\bar{\Delta})=L_\nu(P,\Delta)$ $\forall
\bar{\Delta}\in\ker(V^\star)$, we can assume that
$\Delta\in[\ker(V^\star)]^\bot$. Moreover, $L_\nu(\cdot,\Delta)$
is strictly convex over $\Qc_{n,+}$. Thus, the unique minimum
point of $L_\nu(\cdot,\Delta)$, which is given by annihilating the
first directional derivative of $L_\nu(\cdot,\Delta)$, is \eq
P_\nu(\Delta):=\left(\hat{\Sigma}_C^{-\frac{1}{\nu}}+\frac{1}{\nu}V^\star(\Delta)\right)^{-\nu}.\nn\eeq
Since $P_\nu(\Delta)\in\Qc_{n,+}$, the set of the admissible {\em
Lagrange} multipliers is \eq \Lc_\nu:=\Set{\Delta\in\Qc_n\;|\;
\hat{\Sigma}_C^{-\frac{1}{\nu}}+\frac{1}{\nu}V^\star(\Delta)>0}\cap
[\ker(V^\star)]^\bot\nn\eeq which is an open and bounded set (the
proof is similar to the one of Proposition 5.1 in
\cite{ME_ENHANCEMENT_FERRANTE_2012}). Then, the dual problem is
 \eq \Delta^\circ=\Argmin{\Delta\in\Lc_\nu} J_\nu(\Delta)\nn \eeq where
\eqn \label{funzionale_duale_cov} J_\nu(\Delta)&:=&
-\inf_{P}L_\nu(P,\Delta)\nn\\&=&\frac{\nu}{\nu-1}\tr\left(\hat{\Sigma}_C^{-\frac{1}{\nu}}+\frac{1}{\nu}V^\star(\Delta)\right)^{1-\nu}.\eeqn
Note that
$J_\nu(0)=\frac{\nu}{\nu-1}\tr\left[\hat{\Sigma}_C^{\frac{\nu-1}{\nu}}\right]$.
Accordingly, we can restrict the search of a minimum point to the
set $\Lc^\star:=\Set{\Delta\in\Lbg\;|\; J_\nu(\Delta)\leq
J_\nu(0)}\subset \Lc_\nu$ which is bounded. Following the same
lines in \cite{ME_ENHANCEMENT_FERRANTE_2012}, it is possible to
prove that $J_\nu\in\Cc^\infty(\Lc_\nu)$ is strictly convex on
$\Lc_\nu$ and $\Lim{\Delta\rightarrow
\partial \Lc_\nu} J_\nu(\Delta)=+\infty$ (the limit
diverges because the exponent in (\ref{funzionale_duale_cov}) is
negative). Thus, $\Lc^\star$ is a compact set (i.e. closed and
bounded) and $J_\nu$ admits a minimum point $\Delta^\circ$ over
$\Lc^\star$ by the Weierstrass' Theorem. The uniqueness of
$\Delta^\circ$ follows from the fact that $J_\nu$ is strictly
convex over $\Lc_\nu$. Also in this case, a globally convergent
matricial {\em Newton} algorithm for finding $\Delta^\circ$ may be
employed. Therefore, once we computed $\Delta^\circ$ the solution
to Problem \ref{problema_ancillario} is given by
$P_\nu(\Delta^\circ)$. Finally, the same analysis may be extended
to $\Dc_{\mathrm{KL}}$. In this case,
$P_{KL}(\Delta)=\e{\log(\hat{\Sigma}_C)-V^\star(\Delta)}$.

To sum up, a family of solutions $P_\nu\in\Rgamma\cap \Qc_{n,+}$
to the structured covariance estimation problem has been found. In
this way, we have a complete tool to compute the family of
estimates $\Phi_\nu$ of $\Omega$ starting from a {\em prior}
$\Psi$ and a finite-length data $y_1 \ldots y_N$: We compute
$P_\nu$ from $y_1 \ldots y_N$ and we then find $\Phi_\nu$ starting
from $P_\nu$ and $\Psi$.

\section{Simulation results - Part II}\label{section_sim2}
We consider the bivariate {\em bandpass} random process $y$ of
Section \ref{subsection_simu_multi} and we take into account the
following THREE-like spectral estimation procedure:
\begin{enumerate}
    \item We start from a finite sequence $y_1 \ldots y_N$
    extracted from a realization of the process $y$;
    \item Fix $G(z)$ as in Section \ref{subsection_simu_multi};
    \item Choose a prior spectral density $\Psi\in\Stm$ such that $\Psi^{\frac{1}{\nu}}$ is rational;
    \item Feed the filters bank with the data sequence $y_1\ldots y_N$, collect the output data $x_1\ldots x_N$ and compute  $\hat{\Sigma}_C=\frac{1}{N}\sum_{k=1}^N x_kx_k^T$;
    \item Compute $P_\nu \in \Rgamma\cap\Qc_{n,+}$ by solving Problem
    \ref{problema_ancillario} (with $\Dc_\nu$), then set $\hat{\Sigma}=P_\nu$;
    \item Compute $\Phi_\nu$ by solving Problem \ref{problema_THREE_I} (with $\Dnu$) by means of the
    proposed algorithm with the chosen $\Psi$ and  $\hat{\Sigma}^{-\frac{1}{2}}G(z)$ as filters bank.
\end{enumerate}
As noticed in Section \ref{section_sim}, $\Psi$ represents the {\em
a priori} information on $y$. Accordingly, $\Phi_\nu$ is a
spectral density (with bounded {\em McMillan} degree) which is
consistent with the interpolation constraint in
(\ref{problema_THREE}) and is as close as possible to the {\em a
priori} information, encoded in $\Psi$, according to the
divergence index $\Dnu$.

In the following example, one can consider the prior in Figure
\ref{foto2}. If no {\em a priori} information is given, we set
$\Psi=I$. However, one can get information on $y$ by choosing
$\Psi$ as the constant spectral density equal to the variance of
the given data sequence. In this way the corresponding estimate
will possess at least the zeroth moment similar to the estimated
one by the given data. In what follows, the latter has been
considered. In Figure \ref{foto4}, the obtained estimates with
$N=50$ (i.e. we have considered a short-length data) are depicted.
For the extracted sequence, the estimators for $\nu=2$ and $\nu=3$
appear to perform better than the one for $\nu=1$. More precisely,
the peaks of the estimates are reduced by increasing $\nu$. In
fact, as shown in Section \ref{section_sim}, large values of $\nu$
penalize solutions which are very different from $\Psi$ in narrow
ranges of frequencies. In this case $\Psi$ is constant, thus
solutions with $\nu$ large will be more ``flat'' than the one with
$\nu=1$.

\begin{figure}[htbp]
\begin{center}
\includegraphics[width=8.8cm]{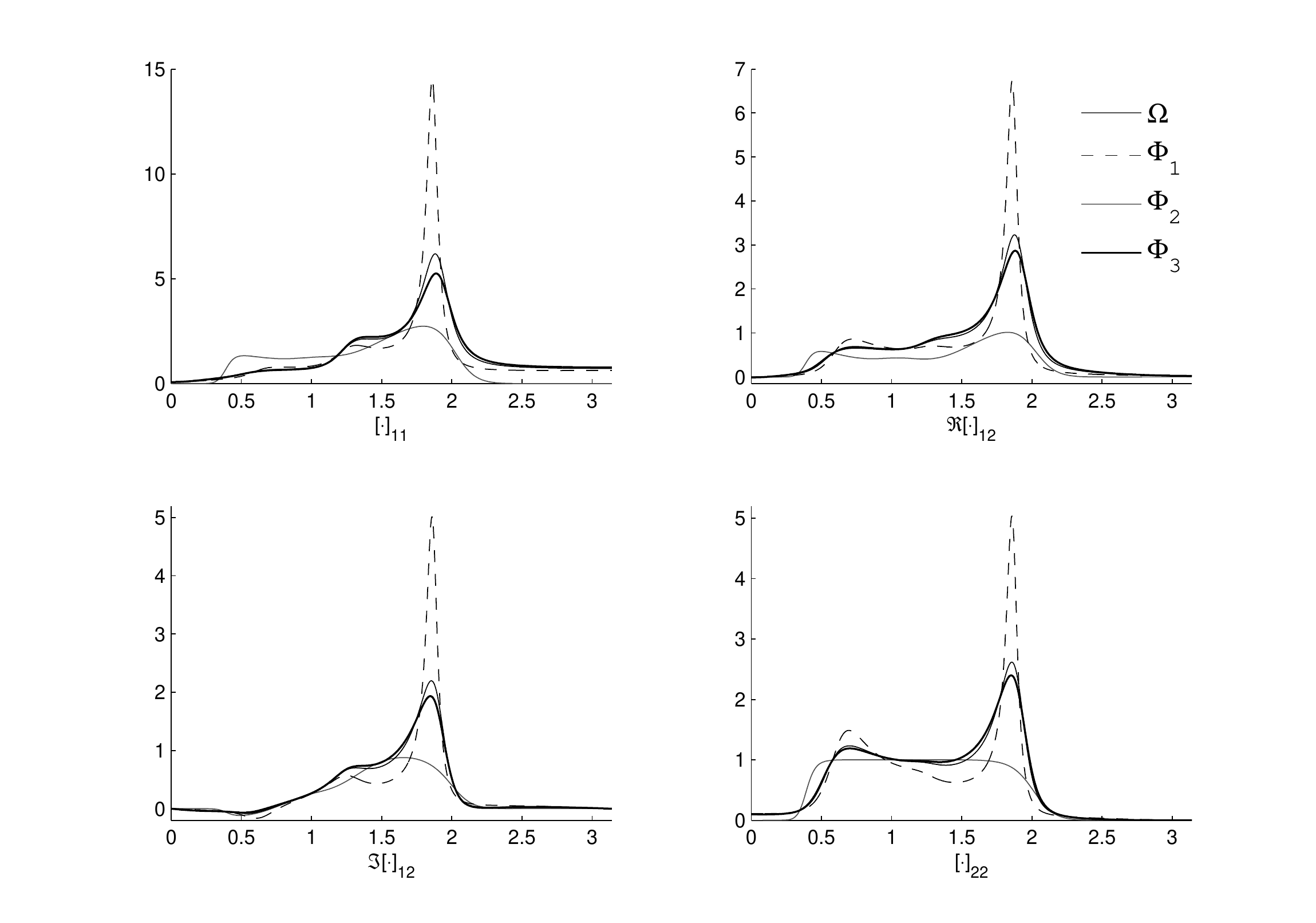}
\end{center}
\vspace{-6mm} \caption{Estimation of the spectral density of a
bivariate {\em bandpass} random process.}\label{foto4}
\end{figure}

In the light of the results found in Section \ref{section_sim} and
here, we can outline the application scenarios for the presented
family of estimators. The estimator with $\nu=1$ is preferable
when the {\em a priori} model for $y$ is rational (i.e. $\Psi$
rational) and a simple model for $y$ is required. On the contrary,
estimators with $\nu$ large are preferable when the model for $y$
must be similar to the {\em a priori} model also in narrow ranges
of frequencies and it must exhibit a ``rich'' dynamic. Increasing
$\nu$, the previous features become more remarked. The limit
case is $\nu\rightarrow\infty$ and the corresponding model is generically non-rational.

\section{Conclusions}
A multivariate Beta divergence family connecting the {\em
Itakura-Saito} distance with the {\em Kullback-Leibler} divergence
has been introduced. The corresponding solutions to the spectrum
approximation problem are rational
 when the parametrization in Theorem \ref{teo_riassuntivo} of the parameter $\beta$ is employed. Such family
also includes the solution corresponding to the {\em
Itakura-Saito} distance. Moreover, the limit of this family tends
to the solution corresponding to the {\em Kullback-Leibler}
divergence. Then, similar results may be found for the structured
covariance estimation problem. Simulations, together with the
potential application scenarios, suggest that the presented family
of estimators provides a relevant tool in multivariate spectral
estimation.

\section*{acknowledgments}
The author wish to thank Prof. A. Ferrante and Prof. M. Pavon for
useful discussions about THREE-like spectral estimation
approaches.

\appendix
\subsection{On the exponentiation of positive definite matrices}
\label{exp_matrice}We collect some technical result concerning the
exponentiation of positive definite matrices to an arbitrary real
number. We start by introducing the differential of the matrix
exponential and the matrix logarithm (see
\cite{RELATIVE_ENTROPY_GEORGIOU_2006}). \prop Given $Y\in\Qc_{n}$,
the differential of $Y\mapsto \e{Y}$ in the direction
$\Delta\in\Qc_n$ is given by the linear map \eq M_Y: \Delta\mapsto
\int_0^1\e{(1-\tau)Y}\Delta\e{\tau Y}\de \tau.\nn\eeq \eprop

\prop \label{dif_logx}Given $Y\in\Qc_{n,+}$, the differential of
$Y\mapsto\log{(Y)}$ in the direction $\Delta\in\Qc_n$ is given by
the linear map \eq N_Y: \Delta\mapsto
\int_0^\infty(Y+tI)^{-1}\Delta(Y+tI)^{-1}\de t.\nn\eeq \eprop  Let us
consider now a positive definite matrix $X\in\Qc_{n,+}$ and a real
number $c$. The exponentiation of $X$ to $c$ may be rewritten in
the following way \eq X^c=e^{c\log X}.\nn\eeq Accordingly, by
applying the chain rule, the differential of $X\mapsto X^c$ in the
direction $\Delta\in\Qc_n$ is given by \eqn &&
\hspace{-0.5cm}M_{c\log X}(cN_X(\Delta))
=c\int_0^1X^{c(1-\tau)}\int_0^\infty (X+tI)^{-1}\nn\\ && \times \Delta
(X+tI)^{-1}\de tX^{c\tau} \de\tau.\nn\eeqn We summarize this result
below. \prop \label{diff_Xc} The differential of $X\mapsto X^c$ in
direction $\Delta\in\Qc_n $ is given by the linear map
\eqn\label{formula_Pxc} && \hspace{-1cm}O_{X,c}:\Delta\mapsto
c\int_0^1X^{c(1-\tau)}\int_0^\infty (X+tI)^{-1}\nn\\ && \times \Delta
(X+tI)^{-1}\de tX^{c\tau} \de\tau. \eeqn\eprop \cor
\label{variazione_Xc} The first variation of $X\mapsto\tr[X^c]$ in
direction $\Delta\in\Qc_n $ is \eq\label{variaz_traccia_Xc}
\delta(\tr[X^c];\Delta)=c\tr(X^{c-1}\Delta). \eeq\ecor \proof
Since $X^{c\tau}$ and $(X+tI)$ commute, we get \eqn
\delta(\tr[X^c];\Delta)&=&\tr(O_{X,c}(\Delta))\nn\\&=&c\tr\left[
X^{c}\int_0^\infty (X+tI)^{-2}\de t\Delta
\right]\nn\\&=&c\tr\left[ X^{c}X^{-1}\Delta \right]=c\tr\left[
X^{c-1}\Delta \right].\nn\eeqn \qed

\subsection{Proofs of Section \ref{section_beta_family}}\label{proofs_SEC_III}

{\em Proof of Proposition \ref{prop_continuity_discrepancy}}: By definition
$X^{1-c}$ and $Y^{c-1}$ are continuous function of real variable
$c$. Thus, the function $\loga{c}(X,Y)$ of real variable $c$ is
continuous in $\Rs\setminus\Set{1}$. It remains to prove that
$\loga{c}$ is continuous in $c=1$. This is equivalent to show that
$\lim_{c\rightarrow 1} \loga{c}(X,Y)=\log(X)-\log(Y)$. Let
$X=U\mathrm{diag}(d_1,\ldots, d_m)U^T$, then \eq
\frac{1}{1-c}(X^{1-c}-I)=U\mathrm{diag}\left(\frac{d_1^{1-c}-1}{1-c},\ldots,\frac{d_m^{1-c}-1}{1-c}\right)U^T.
\nn \eeq Taking the limit for $c\rightarrow 1$, we get \eqn
\label{limite_logaritmo}&&\lim_{c\rightarrow 1}
\frac{1}{1-c}(X^{1-c}-I)\nn \\
&&\hspace{+0.4cm}=U\mathrm{diag}\left(\lim_{c\rightarrow
1}\frac{d_1^{1-c}-1}{1-c},\ldots,\lim_{c\rightarrow
1}\frac{d_m^{1-c}-1}{1-c}\right)U^T \nn\\ &&\hspace{+0.4cm}
=U\mathrm{diag}\left(\log(d_1),\ldots,\log(d_m)\right)U^T=\log(X).
\eeqn
Accordingly,  \eqn \label{limite_logaritmo_discrepancy}&&\lim_{c\rightarrow 1}\loga{c}(X,Y)\nn\\
&&\hspace{+0.4cm}=\lim_{c\rightarrow 1}
\left(\frac{1}{1-c}(X^{1-c}-I)-\frac{1}{1-c}(Y^{1-c}-I)\right)Y^{c-1}\nn\\&&\hspace{+0.4cm}=
\lim_{c\rightarrow
1}\left(\frac{1}{1-c}(X^{1-c}-I)\right)-\lim_{c\rightarrow
1}\left(\frac{1}{1-c}(Y^{1-c}-I)\right)\nn\\
&&\hspace{+0.4cm}=\log(X)-\log(Y)\eeqn which proves that
$\loga{c}$ is continuous in $c=1$. Concerning the last statement,
it is straightforward that $X=Y$ implies $\log_c(X,Y)=0$. On the
contrary, $\log_{c}(X,Y)=0$, with $c\neq 1$, implies
$X^{1-c}Y^{c-1}=I$ which is equivalent to $X^{1-c}=Y^{1-c}$, since
$X,Y\in\Qc_{m,+}$. Thus, $X=Y$.  We get the same
conclusion for $c=1$ by using similar argumentations.\qed\\

 {\em Proof of Proposition \ref{prop_limiti_beta_multi}}: Since $\Phi$ and $\Psi$ belong to $\Stm$, i.e.
$\Phi$ and $\Psi$ are coercive and bounded, it is possible to show
by standard argumentations that the integrand function of
(\ref{def_beta_multi}) uniformly converges on $\Ts$ for
$\beta\rightarrow 0$ and $\beta\rightarrow 1$. Hence, it is
allowed to pass the limits, for $\beta\rightarrow 0$ and
$\beta\rightarrow 1$, under the integral sign. Taking into account
the first limit, we get \eqn && \hspace{-0.4cm}\lim_{\beta\rightarrow
0}\Db(\Phi\|\Psi)\nn\\ &&\hspace{-0.1cm}=\lim_{\beta\rightarrow
0}\int\tr\left[
\frac{1}{\beta-1}(\Phi^\beta-\Phi\Psi^{\beta-1})-\frac{1}{\beta}(\Phi^{\beta}-\Psi^\beta)\right]\nn\\
&&\hspace{-0.1cm}=\int\tr\left[
-I+\Phi\Psi^{-1}-\lim_{\beta\rightarrow
0}\frac{1}{\beta}\left((\Phi^{\beta}-I)-(\Psi^\beta-I)\right)\right]\nn\\
&&\hspace{-0.1cm}=\int\tr\left[-I+\Phi\Psi^{-1}-
\log(\Phi)+\log(\Psi)\right]\nn\\&&\hspace{-0.1cm}=\Dis(\Phi\|\Psi)\nn\eeqn
where we exploited (\ref{limite_logaritmo}). For the second limit,
we obtain
 \eqn && \lim_{\beta\rightarrow 1}\Db(\Phi\|\Psi)\nn\\&&\hspace{+0.2cm}=\lim_{\beta\rightarrow 1}\left(-\frac{1}{\beta}\int\tr\left[\Phi^{\beta}
 \loga{\frac{1}{\beta}}\left(\Psi^{\beta},\Phi^{\beta}\right)
+\Phi^{\beta}-\Psi^{\beta}\right]\right)\nn\\
&&\hspace{+0.2cm}=-\int\tr\left[\Phi\lim_{\beta\rightarrow
1}\loga{\frac{1}{\beta}}\left(\Psi^{\beta},\Phi^{\beta}\right)
+\Phi-\Psi\right]\nn\\ &&
\hspace{+0.2cm}=-\int\tr\left[\Phi\lim_{\beta\rightarrow
1}\loga{2-\beta}\left(\Psi,\Phi\right) +\Phi-\Psi\right]\nn\\
&&\hspace{+0.2cm}=\int\tr\left[\Phi\left(\log(\Phi)-\log(\Psi)\right)
+\Psi-\Phi\right]\nn\\&&\hspace{+0.2cm}=\Dkl(\Phi\|\Psi)\nn\eeqn where we exploited (\ref{limite_logaritmo_discrepancy}).\qed\\

{\em Proof of Proposition \ref{proposizione_beta_div_strict_convex}}: The proof will be divided in
the following three cases: $0<\beta<1$, $\beta=1$ and
$\beta=0$.\\{\em $\bullet$ Case $0<\beta<1$: Point 1.} The first
variation of $\Db(\Phi\|\Psi)$, with respect to $\Phi$, in
direction $\delta\Phi\in \Lm$ is \eqn \label{prima_var_Db}&&
\delta(\Db(\Phi\|\Psi);\delta \Phi)\nn\\ &&\hspace{0.5cm}=\frac{1}{\beta-1}\int_0^{2\pi}
\tr\left[(\Phi^{\beta-1}-\Psi^{\beta-1})\delta
\Phi\right]\frac{\de\vartheta}{2\pi}\eeqn where we exploited (\ref{variaz_traccia_Xc}). The second variation in
direction $\delta\Phi$ is \eqn  &&\delta^2(\Db(\Phi\|\Psi);\delta
\Phi)=\int_0^{2\pi}
\tr\left[O_{\Phi,\beta-1}(\delta\Phi)\delta\Phi\right]\frac{\de\vartheta}{2\pi}\nn\\
&& \hspace{0.5cm}=\int_0^{2\pi}\tr\left[\int_0^1\Phi^{(\beta-1)(1-\tau)}\int_0^\infty
(\Phi+tI)^{-1}\delta\Phi \right. \nn\\ &&\hspace{+0.8cm}\times \left.(\Phi+tI)^{-1}\de t\Phi^{(\beta-1)\tau}\de \tau \delta\Phi\right]\frac{\de\vartheta}{2\pi}\nn\\
&&\hspace{0.5cm}=\int_0^{2\pi}\int_0^1\int_0^\infty
\tr\left[\Phi^{(\beta-1)(1-\tau)}(\Phi+tI)^{-1}\delta\Phi \right.\nn\\
&&\hspace{+0.8cm}\left.\times
(\Phi+tI)^{-1}\Phi^{(\beta-1)\tau}\delta\Phi\right]\de t\de
\tau\frac{\de\vartheta}{2\pi} \nn\eeqn  where $O_{X,c}$ is defined in (\ref{formula_Pxc}).
By the cyclic property of
the trace and since $\Phi^{(\beta-1)\tau}$ and $(\Phi+tI)^{-1}$
commute, we get \eq
\label{integrale_sec_variazizone_convex_dist}\delta^2(\Db(\Phi\|\Psi);\delta
\Phi)=\int_0^{2\pi}\int_0^1\int_0^\infty
f_{t,\tau}(\Phi,\delta\Phi)\de t \de \tau \frac{\de
\vartheta}{2\pi}\eeq where \eqn && \hspace{-0.8cm}
f_{t,\tau}(X,\Delta)=\tr\left[X^{\frac{(\beta-1)\tau}{2}}(X+tI)^{-\frac{1}{2}}\Delta
(X+tI)^{-\frac{1}{2}}\right.\nn\\ && \hspace{-0.4cm}\left.\times
X^{(\beta-1)(1-\tau)} (X+tI)^{-\frac{1}{2}} \Delta
(X+tI)^{-\frac{1}{2}}X^{\frac{(\beta-1)\tau}{2}}\right]\nn\eeqn with
$X\in\Qc_{m,+}$, $\Delta\in\Qc_m$, $t\in[0,\infty)$ and
$\tau\in[0,1]$. Thus, $f_{t,\tau}(X,\Delta)\geq 0$ and
$f_{t,\tau}(X,\Delta)= 0$ if and only if $\Delta= 0$. We conclude
that integral (\ref{integrale_sec_variazizone_convex_dist}), i.e.
the second variation of $\Db(\cdot\|\Psi)$, is positive for
$\delta\Phi\neq 0$. Accordingly, $\Db(\cdot\|\Psi)$
is strictly convex over the convex set $\Stm$.\\
{\em Point 2.} As a consequence of the previous statement, the minimum point
is unique and it is given by annihilating (\ref{prima_var_Db}) for
each $\delta\Phi\in\Lm$. Since
$\Phi^{\beta-1}-\Psi^{\beta-1}\in\Lm$, it follows that the minimum
point satisfies the condition $\Phi^{\beta-1}=\Psi^{\beta-1}$.
Accordingly, $\Phi=\Psi$. Finally it is sufficient to observe that
$\Db(\Psi\|\Psi)=0$.\\ {\em $\bullet$ Case $\beta=1$}: Firstly,
given $X\in\Qc_{m,+}$ and $\delta X\in\Qc_m$ we have \eqn
\label{formule_var_log} \delta(\tr[\log(X)];\delta
X)&=&\tr[X^{-1}\delta X]\nn\\\delta(\tr[X\log(X)];\delta
X)&=&\tr[(\log(X)+I)\delta X]\eeqn (it is sufficient to apply
Proposition \ref{dif_logx} in Appendix \ref{exp_matrice}). For
$\beta=1$ we get the {\em Kullback-Leibler} divergence in
(\ref{dist_KL}). Taking into account (\ref{formule_var_log}), its
first and second variations with respect to $\Phi$ in direction
$\delta\Phi\in L_\infty^{m\times m}(\Ts)$ are, respectively, \eqn
\delta(\Dkl(\Phi\|\Psi);\delta \Phi)&=&\int
\tr\left[(\log(\Phi)-\log(\Psi))\delta \Phi\right]\nn\\
\delta^2(\Dkl(\Phi\|\Psi);\delta \Phi)&=&\int \tr\left[\delta
\Phi\Phi^{-1}\delta \Phi\right].\nn\eeqn Since the second variation
is non negative and equal to zero if and only if $\delta
\Phi\equiv 0$, $\Dkl(\cdot\|\Psi)$ is strictly convex over $\Stm$
and the (unique) minimum point is given by annihilating the first
variation which leads to condition $\log (\Phi)=\log(\Psi)$. Thus,
$\Phi=\Psi$ and $\Dkl(\Psi\|\Psi)=0$.\\ {\em $\bullet$ Case
$\beta=0$}: In this case we have the {\em Itakura-Saito} distance.
Using similar argumentations used for the case $\beta=1$, we
get the
statement.\qed\\

\end{document}